\documentclass[11pt, a4paper]{amsart}
\usepackage{amsmath, amsthm, amssymb}
\usepackage{txfonts, mathrsfs}
\usepackage[usenames]{color}
\usepackage{psfrag}
\usepackage{graphicx}

\numberwithin{equation}{section}

\newtheorem{thm}{Theorem}[section]
\newcommand{\bt}{\begin{thm}}
\newcommand{\et}{\end{thm}}

\newtheorem{cor}[thm]{Corollary}
\newcommand{\bc}{\begin{cor}}
\newcommand{\ec}{\end{cor}}

\newtheorem{lem}[thm]{Lemma}
\newcommand{\bl}{\begin{lem}}
\newcommand{\el}{\end{lem}}

\newtheorem{prop}[thm]{Proposition}
\newcommand{\bp}{\begin{prop}}
\newcommand{\ep}{\end{prop}}

\newtheorem{defn}[thm]{Definition}
\newcommand{\bd}{\begin{defn}}
\newcommand{\ed}{\end{defn}}

\newtheorem{rmrk}[thm]{Remark}
\newcommand{\br}{\begin{rmrk}}
\newcommand{\er}{\end{rmrk}}

\newtheorem{quest}[thm]{Question}
\newcommand{\bq}{\begin{quest}}
\newcommand{\eq}{\end{quest}}

\newcommand{\C}{\mathbb{C}}
\newcommand{\N}{\mathbb{N}}

\newcommand{\SL}{\operatorname{SL}}
\newcommand{\Hyp}{\mathbb{H}}

\newdimen\vintkern\vintkern12pt
\def\vint{-\kern-\vintkern\int}

\newcommand{\hm}{{\mathcal H}}

\newcommand{\trace}{\operatorname{tr}}
\newcommand{\length}{\ell}
\newcommand{\Area}{\operatorname{Area}}

\newcommand{\md}{\operatorname{md}}

\newcommand{\sys}{\operatorname{sys}}

\newcommand{\jac}{{\mathbf J}}
\newcommand{\ap}{\operatorname{ap}}
\newcommand{\apmd}{\ap\md}

\newcommand{\R}{\mathbb{R}}

\begin{document}
\pagebreak
\bibliographystyle{plain}


\title{Area minimizing surfaces of bounded genus in metric spaces}


\author{Martin Fitzi}

\address
  {Department of Mathematics\\ University of Fribourg\\ Chemin du Mus\'ee 23\\ 1700 Fribourg, Switzerland}
\email{martin.fitzi@unifr.ch}

\author{Stefan Wenger}

\address
  {Department of Mathematics\\ University of Fribourg\\ Chemin du Mus\'ee 23\\ 1700 Fribourg, Switzerland}
\email{stefan.wenger@unifr.ch}

\date{\today}

\thanks{Research supported by Swiss National Science Foundation Grants 165848 and 182423}

\begin{abstract}
 The Plateau-Douglas problem asks to find an area minimizing surface of fixed or bounded genus spanning a given finite collection of Jordan curves in Euclidean space. In the present paper we solve this problem in the setting of proper metric spaces admitting a local quadratic isoperimetric inequality for curves. We moreover obtain continuity up to the boundary and interior H\"older regularity of solutions. Our results generalize corresponding results of Jost and Tomi-Tromba from the setting of Riemannian manifolds to that of proper metric spaces with a local quadratic isoperimetric inequality. The special case of a disc-type surface spanning a single Jordan curve corresponds to the classical problem of Plateau, in proper metric spaces recently solved by Lytchak and the second author.
\end{abstract}

\maketitle

\renewcommand{\theequation}{\arabic{section}.\arabic{equation}}
\pagenumbering{arabic}

\section{Introduction and statement of results}\label{sec:Intro}

\subsection{Introduction} The classical problem of Plateau is about the existence of an area minimizing disc spanning a given Jordan curve in Euclidean $\R^n$. This problem has a long and rich history and admits many variations. The first rigorous solutions go back to Douglas \cite{Dou31} and Rad\'o \cite{Rad30}. Morrey \cite{Mor48} extended the solution to a large class of Riemannian manifolds. Since then many variants of the problem have been studied. Douglas considered in \cite{Dou39} the more general problem of finding an area minimizing surface of fixed genus spanning a finite collection of Jordan curves in $\R^n$. This problem is nowadays called the Douglas or Plateau-Douglas problem. Solutions to this problem were proposed by Douglas \cite{Dou39}, Shiffman \cite{Shi39}, and Courant \cite{Cou40}. Tomi-Tromba \cite{TT88} and Jost \cite{Jos85} solved the Plateau-Douglas problem for a large class of Riemannian manifolds called homogeneously regular. We refer to \cite{DHT10} for a more detailed history of these problems.

Recently, Lytchak and the second author solved the classical problem of Plateau in any proper metric space \cite{LW15-Plateau}. They moreover proved continuity up to the boundary and interior H\"older regularity of area minimizing discs when the underlying metric space admits a local quadratic isoperimetric inequality for curves. See \cite{Nik79}, \cite{MZ10}, \cite{OvdM14} for some earlier results for special classes of metric spaces.

In the present paper we solve the Plateau-Douglas problem in the setting of metric spaces. Firstly, we prove existence of area minimizing surfaces of fixed topological type spanning a finite collection of rectifiable Jordan curves in proper metric spaces admitting a local quadratic isoperimetric inequality. For this we assume the so-called Douglas condition, also used by Douglas \cite{Dou39} and Jost \cite{Jos85}. Secondly, we find area minimizing surfaces (possibly disconnected) of bounded genus without invoking any Douglas type condition. Moreover, we obtain continuity up to the boundary and interior H\"older regularity of solutions. Finally, we discuss the existence of energy minimizers under Courant's condition of cohesion used in \cite{Shi39}, \cite{Cou40}, \cite{TT88}.

Plateau's problem can be considered for other notions of surfaces. For example, in the class of integral currents - generalized surfaces of arbitrary genus - Plateau's problem admits a solution in any compact metric space by \cite{AK00}. Regularity of area minimizing integral currents is however only known in the setting of Riemannian manifolds.

\subsection{Statement of results} In this paper, a smooth surface refers to a smooth orientable two-dimensional manifold, possibly disconnected and with boundary. Let $X$ be a complete metric space and $M$ a smooth compact surface.
For $q>1$ we let $W^{1, q}(M, X)$ be the space of $q$--Sobolev maps from $M$ to $X$ in the sense of Reshetnyak \cite{Res97}, \cite{Res06}. There exist several different but equivalent definitions of Sobolev maps from Euclidean, Riemannian or even more general domains into a metric space. For references as well as for definitions of the following notions see Section~\ref{sec:prelims} below. 
The parametrized Hausdorff area of a map $u\in W^{1,2}(M, X)$ is denoted $\Area(u)$. If $u$ is Lipschitz and injective then $\Area(u)$ is simply the Hausdorff $2$--measure of the image of $u$. When $\partial M$ is non-empty then $u\in W^{1, q}(M, X)$ has a trace $\trace(u)$, which is a map defined almost everywhere on $\partial M$. If $u$ has a representative $\bar{u}$ which is continuous up to the boundary then $\trace(u)$ agrees with the restriction $\bar{u}|_{\partial M}$ almost everywhere on $\partial M$.

\bd\label{def:loc-quad-isop} 
A complete metric space $X$ is said to admit a local quadratic isoperimetric inequality if there exist $C, l_0>0$ such that every Lipschitz curve $c\colon S^1\to X$ of length $\length(c)\leq l_0$ is the trace of a Sobolev map $u\in W^{1,2}(D, X)$ with $$\Area(u)\leq C\cdot \length^2(c),$$ where $D$ denotes the Euclidean unit disc.
\ed

 Spaces admitting a local quadratic isoperimetric inequality include homogeneously regular Riemannian manifolds, compact Lipschitz manifolds and thus all compact Finsler manifolds; moreover all Banach spaces, complete ${\rm CAT}(\kappa)$ spaces for all $\kappa\in\R$, compact Alexandrov spaces, and many more (see \cite{LW15-Plateau}).

Let $M$ be a smooth compact surface with $k\geq 1$ boundary components. Given a disjoint union $\Gamma$ of $k$ Jordan curves in $X$ we let $\Lambda(M,\Gamma,X)$ be the possibly empty family of Sobolev maps $u\in W^{1,2}(M, X)$ such that $\trace(u)$ has a continuous representative which weakly monotonically parametrizes $\Gamma$. A weakly monotone parametrization of $\Gamma$ is a map from $\partial M$ to $\Gamma$ which is the uniform limit of homeomorphisms $\varphi_j\colon \partial M\to \Gamma$. We set $$a(M,\Gamma, X):= \inf\{\Area(u): u\in\Lambda(M, \Gamma, X)\},$$ where the infimum of the empty set is infinite. The value $a(M, \Gamma, X)$ does not change when $M$ is replaced by a surface diffeomorphic to $M$.
If $p\geq 0$ and $M$ is the smooth compact and {\it connected} surface with $k$ boundary components and of genus $p$ (such $M$ is unique up to diffeomorphism) then we abbreviate $$a_p(\Gamma, X):= a(M, \Gamma, X).$$

In our first result we will impose the so-called Douglas condition introduced in \cite{Dou39}. Define $$a_p^*(\Gamma, X):=\min a(M^*,\Gamma, X),$$ where the minimum runs over all compact surfaces $M^*$ with $k$ boundary components having the following property. Either $M^*$ is connected and has genus at most $p-1$ or $M^*$ has exactly two connected components, each of which has non-empty boundary, and the genus of $M^*$ is $p$. Notice that there are only finitely many such $M^*$ up to diffeomorphism. It can be shown that $a_p(\Gamma, X) \leq a_p^*(\Gamma, X)$ whenever the space $X$ is rectifiably connected.
We say that the Douglas condition holds for $p\geq 0$ and $\Gamma$ if we have the strict inequality
 \begin{equation}\label{eq:Douglas-condition}
  a_p(\Gamma, X)<a_p^*(\Gamma, X).
 \end{equation}
This condition thus requires that it is strictly easier to fill $\Gamma$ with a connected surface of genus $p$ than it is to fill $\Gamma$ with a surface of lower topological type.

\bt\label{thm:Douglas-Plateau-intro}
 Let $X$ be a proper metric space admitting a local quadratic isoperimetric inequality, and let $\Gamma$ be the disjoint union of $k\geq 1$ rectifiable Jordan curves in $X$. Let $M$ be a smooth compact and connected surface with $k$ boundary components and of genus $p\geq 0$. If the Douglas condition \eqref{eq:Douglas-condition} holds then there exist $u\in\Lambda(M, \Gamma, X)$ and a Riemannian metric $g$ on $M$ such that $$\Area(u) = a_p(\Gamma, X)$$ and $u$ is infinitesimally isotropic with respect to $g$.
\et

Recall that a metric space is proper if all its closed bounded subsets are compact. The precise meaning of {\it infinitesimally isotropic} will be given in Section~\ref{sec:energy-min-isotropic}. It provides a substitute for a conformal parametrization in the setting of (non-Euclidean) metric spaces. If $u$ is infinitesimally isotropic then it is in particular $\sqrt{2}$--quasiconformal: it maps infinitesimal balls in $(M,g)$ to `ellipses' of eccentricity at most $\sqrt{2}$. If $X$ is a Riemannian manifold or, more generally, a metric space with property (ET) in the sense of \cite{LW15-Plateau} then $\sqrt{2}$ can be replaced by $1$. 

The Riemannian metric $g$ in Theorem~\ref{thm:Douglas-Plateau-intro} can be chosen in such a way that $(M,g)$ has constant curvature $-1$, $0$, or $1$ and that $\partial M$ is geodesic. When $k=1$ and $p=0$ then the theorem asserts the existence of an area minimizing disc which was treated in \cite{LW15-Plateau}. Theorem~\ref{thm:Douglas-Plateau-intro} generalizes the corresponding results in \cite{Jos85} and \cite{TT88} from the setting of homogeneously regular Riemannian manifolds to that of proper metric spaces with a local quadratic isoperimetric inequality.

Our next result is concerned with the existence of (possibly disconnected) minimal surfaces of bounded genus in the same setting. Let $\hat{\mathcal{M}}(k, p)$ be the family of smooth compact surfaces $\hat{M}$ with $k$ boundary components and of genus at most $p$ such that each connected component of $\hat{M}$ has non-empty boundary. We set $$\hat{a}_p(\Gamma, X):= \min\left\{ a(\hat{M}, \Gamma, X): \hat{M}\in \hat{\mathcal{M}}(k, p)\right\}$$ and notice that $\hat{a}_p(\Gamma, X) \leq a_p(\Gamma, X)$, with equality for example if $X$ is rectifiably connected. 

\bt\label{thm:area-min-bounded-genus}
 Let $X$ be a proper metric space admitting a local quadratic isoperimetric inequality. Let $\Gamma$ be the disjoint union of $k\geq 1$ rectifiable Jordan curves in $X$ and let $p\geq 0$. If $\hat{a}_p(\Gamma, X)<\infty$ then there exist $\hat{M}\in \hat{\mathcal{M}}(k,p)$ and $u\in \Lambda(\hat{M}, \Gamma, X)$ as well as a Riemannian metric $g$ on $\hat{M}$ such that $$\Area(u) = \hat{a}_p(\Gamma, X)$$ and $u$ is infinitesimally isotropic with respect to $g$.
\et

The Riemannian metric $g$ on $\hat{M}$ can be chosen in such a way that $\partial \hat{M}$ is geodesic and each connected component of $(\hat{M},g)$ has constant curvature $-1$, $0$, or $1$. Theorem~\ref{thm:area-min-bounded-genus} will be derived from Theorem~\ref{thm:Douglas-Plateau-intro}.

Finally, we turn to the regularity of solutions to the Plateau-Douglas problem. Exactly as in the case of area minimizing discs \cite{LW15-Plateau}, we obtain local H\"older continuity in the interior and continuity up to the boundary. Let $M$ be a smooth compact and possibly disconnected surface with $k\geq 1$ boundary components. 

\bt\label{thm:regularity-sol-Douglas-Plateau}
 Let $X$ be a complete metric space admitting a local quadratic isoperimetric inequality with isoperimetric constant $C$. Let $\Gamma$ be the disjoint union of $k$ Jordan curves in $X$. If $u$ is an area minimizer in $\Lambda(M, \Gamma, X)$ and $u$ is infinitesimally isotropic with respect to some Riemannian metric on $M$ then:
 \begin{enumerate}
  \item[(i)] There exists $q>2$ such that $u\in W^{1,q}_{\rm loc}(M\setminus \partial M, X)$. In particular, $u$ has a representative $\bar{u}$ which is continuous in $M\setminus \partial M$ and satisfies Lusin's property $(N)$.
  \item[(ii)] The representative $\bar{u}$ is locally $\alpha$--H\"older continuous in $M\setminus \partial M$ with $\alpha=(8\pi C)^{-1}$ and extends continuously to $\partial M$.
  \item[(iii)] If every Jordan curve in $\Gamma$ is a chord-arc curve then $\bar{u}$ is H\"older continuous on all of $M$.
  \end{enumerate}
\et

This theorem easily follows from the corresponding regularity results for area minimizing surfaces of disc-type established in \cite{LW15-Plateau}. 

We formulated our results above only for the parametrized Hausdorff area. They moreover hold for the paramatrized area coming from any definition of volume inducing quasi-convex $2$--volume densities. We refer to \cite{LW15-Plateau} for the terminology and for examples from convex geometry.

In our first theorem we used the Douglas condition \eqref{eq:Douglas-condition}. A different condition, called {\it condition of cohesion}, was used by Shiffman \cite{Shi39}, Courant \cite{Cou40}, Tomi-Tromba \cite{TT88}. In Theorem~\ref{thm:cond-cohesion-energy-min} we will show that if there exists an energy minimizing sequence satisfying the condition of cohesion then one can find an energy minimizer in $\Lambda(M,\Gamma, X)$, even when $X$ does not admit an isoperimetric inequality. This generalizes the corresponding results in \cite{Shi39}, \cite{Cou40} and \cite{TT88} to the setting of proper metric spaces. Notice that such energy minimizers need not be minimizers for the parametrized Hausdorff area in the generality of non-Euclidean metric spaces. See Section~\ref{sec:cohesion} for a discussion of this and the existence of area minimizers under the condition of cohesion.

\subsection{Elements of proof} We briefly present some of the steps in the proof of Theorem~\ref{thm:Douglas-Plateau-intro}, which combines methods and ideas from \cite{LW15-Plateau} and \cite{Jos85}. We focus on the case that $M$ has strictly negative Euler characteristic so that $M$ admits a hyperbolic metric $g$. We denote by $E_+^2(u, g)$ the (Reshetnyak) energy  of a map $u\in W^{1,2}(M, X)$ with respect to $g$, see Section~\ref{sec:Sobolev}.

The first ingredient in the proof is Theorem~\ref{thm:energy-min-isotropic} which shows that if $E_+^2(u,g)$ is minimal in the sense that $$E_+^2(u,g)\leq E_+^2(u\circ\varphi, g')$$ for all hyperbolic metrics $g'$ and biLipschitz homeomorphisms $\varphi\colon (M, g')\to (M, g)$ then $u$ is infinitesimally isotropic with respect to $g$. This generalizes the corresponding result proved in \cite{LW15-Plateau} for Sobolev maps defined on the Euclidean unit disc.

We then show in Proposition~\ref{prop:bound-systole} that for all $K,\eta>0$ there exists $\varepsilon>0$ with the following property. If $u\in\Lambda(M,\Gamma, X)$ satisfies 
\begin{equation}\label{eq:intro-proof-mainthm}
\Area(u)\leq a_p^*(\Gamma, X)-\eta
\end{equation}
 and if $g$ is a hyperbolic metric on $M$ such that $E_+^2(u,g)\leq K$ then the relative systole of $(M,g)$ cannot be smaller than $\varepsilon$. For the definition of the relative systole see Section~\ref{sec:rel-sys}. The proof of the proposition relies on the collar lemma from hyperbolic geometry, a Fubini type argument and the local quadratic isoperimetric inequality. Similarly, we prove that for all $K,\eta>0$ the family of continuous representatives of traces of maps $u\in\Lambda(M, \Gamma, X)$ satisfying \eqref{eq:intro-proof-mainthm} and $E_+^2(u, g)\leq K$ is equi-continuous, see Proposition~\ref{prop:equi-cont-traces}. 

Combining Theorem~\ref{thm:energy-min-isotropic} and Propositions~\ref{prop:equi-cont-traces} and \ref{prop:bound-systole} with the Mumford compactness theorem from hyperbolic geometry, the Rellich compactness theorem for Sobolev maps and lower semi-continuity of area and energy  we then establish the following claim. For every $u\in\Lambda(M, \Gamma, X)$ with $\Area(u)<a_p^*(\Gamma, X)$ there exist $v\in \Lambda(M, \Gamma, X)$ and a hyperbolic metric $g$ on $M$ such that $\Area(v)\leq \Area(u)$ and $v$ is infinitesimally isotropic with respect to $g$. This implies in particular that $E_+^2(v, g)$ is bounded by a fixed constant multiple of $\Area(v)$. The map $v$ is obtained by minimizing the energy $E_+^2(u',g')$ over all $u'\in \Lambda(M, \Gamma, X)$ with $\Area(u')\leq \Area(u)$ and all hyperbolic metrics $g'$ and using the results mentioned before the claim. 

Finally, the claim above allows us to consider an area minimizing sequence $(v_n)$ in $\Lambda(M, \Gamma, X)$ together with a sequence $(g_n)$ of hyperbolic metrics on $M$ with the property that the energies $E_+^2(v_n, g_n)$ are uniformly bounded. Arguments similar to that in the proof of the claim then yield the existence of an area minimizer in $\Lambda(M, \Gamma, X)$ which, by the same claim, can be assumed to be infinitesimally isotropic with respect to some hyperbolic metric $g$.

\bigskip

{\bf Acknowledgements:} We thank Patrick Ghanaat, Hugo Parlier and Teri Soultanis for comments and discussions. The second author would moreover like to thank Alexander Lytchak for many fruitful collaborations and inspiring discussions over the past years.

\section{Preliminaries}\label{sec:prelims}

\subsection{Basic notation}

The Euclidean norm of a vector $v\in\R^n$ is denoted $|v|$, the open unit disc in Euclidean $\R^2$ by $$D:= \{z\in\R^2: |z|<1\},$$ and its closure by $\overline{D}$.
The genus $p$ of a compact and connected surface $M$ with $k\geq 0$ boundary components is related to the Euler characteristic $\chi(M)$ of $M$ by the formula $$\chi(M) = 2 - 2p - k.$$ We define the genus of a compact and disconnected surface (possibly with boundary) as the sum of the genera of its connected components.

Let $(X, d)$ be a metric space. The open and closed balls in $X$ centered at $x$ and of radius $r$ are denoted by $B(x,r)$ and $\bar{B}(x,r)$, respectively. 
The $s$--dimensional Hausdorff measure of a subset $A\subset X$ is denoted $\hm_{X}^s(A)$. We choose the normalizing constant in such a way that $\hm_X^n$ coincides with Lebesgue measure when $X$ is Euclidean $\R^n$. In particular, if $(M,g)$ is a Riemannian manifold of dimension $n$ then the Hausdorff $n$--measure $\hm_g^n:=\hm_{(M,g)}^n$ on $(M,g)$ coincides with the Riemannian volume. We write $|A|$ for the Lebesgue measure of a subset $A\subset\R^n$.

\subsection{Sobolev maps with values in metric spaces}\label{sec:Sobolev}

There exist several equivalent theories of Sobolev maps from a Euclidean, Riemannian or even more general domain into a complete metric space, see for example \cite{Amb90}, \cite{KS93}, \cite{Haj96}, \cite{Res97}, \cite{Jos97}, \cite{HKST15}. In what follows we briefly recall the definition introduced by Reshetnyak \cite{Res97}, \cite{Res06}. We furthermore recall the notions of approximate metric differentiability and of parametrized area of a Sobolev map from \cite{LW15-Plateau}.

Let $(X, d)$ be a complete metric space and let $M$ be a smooth compact $n$--dimen\-sional manifold, possibly with non-empty boundary. We will actually only need the cases $n=1$ and $n=2$.
Fix a Riemannian metric $g$ on $M$ and let  $\Omega\subset M$ be an open set and $q>1$. We denote by $L^q(\Omega, X)$ the collection of measurable and essentially separably valued maps $u\colon \Omega\to X$ with the following property. For some and thus every $x\in X$ the function $$u_x(z):= d(x, u(z))$$ belongs to the space $L^q(\Omega)$ of $q$--integrable functions on $\Omega$. A sequence $(u_k)\subset L^q(\Omega, X)$ is said to converge in $L^q(\Omega, X)$ to a map $u\in L^q(\Omega, X)$ if $$\int_\Omega d^q(u_k(z), u(z))\,d\hm_g^n(z)\to 0$$ as $k$ tends to infinity.

\bd
 A map $u\in L^q(\Omega, X)$ belongs to the Sobolev space $W^{1,q}(\Omega, X)$ if
 \begin{enumerate}
  \item[(i)] for every $x\in X$ the function $u_x$ belongs to the Sobolev space $W^{1,q}(\Omega\setminus \partial M)$ of real-valued functions, and
  \item[(ii)] there exists $h\in L^q(\Omega)$ such that for all $x\in X$ we have $|\nabla u_x|\leq h$ almost everywhere on $\Omega$.
 \end{enumerate}
\ed

In the above, $|\nabla u_x|$ is the length of the weak gradient $\nabla u_x$ of $u_x$ with respect to the metric $g$. The space $W^{1,q}_{\rm loc}(\Omega, X)$ is defined in an analogous way.

Let $V\subset \R^n$ be an open set and $z\in V$. A map $v\colon V\to X$ is said to be approximately metrically differentiable at $z$ if there exists a necessarily unique semi-norm $s$ on $\R^n$ such that $$\ap\lim_{z'\to z} \frac{d(v(z'), v(z)) - s(z'-z)}{|z'-z|} = 0,$$ where $\ap\lim$ denotes the approximate limit, see for example \cite{EG92}. If such a semi-norm exists, it is called the approximate metric derivative of $v$ at $z$ and denoted $\apmd v_z$.
Let $\varphi\colon W\to V$ be a diffeomorphism, where $W\subset\R^n$ is some open set, and let $w\in W$. If $v\colon V\to X$ is approximately metrically differentiable at $\varphi(w)$ then the composition $v\circ\varphi$ is approximately metrically differentiable at $w$ with $$\apmd (v\circ\varphi)_w = \apmd v_{\varphi(w)} \circ d\varphi_w.$$

Together with \cite[Proposition 4.3]{LW15-Plateau} this implies that if $u\in W^{1,q}(\Omega, X)$ then for almost every $z\in \Omega$ the composition $u\circ\psi^{-1}$ is approximately metrically differentiable at $\psi(z)$ for some and thus any chart $(U, \psi)$ around $z$. Moreover, the semi-norm on $T_zM$ defined by $$\apmd u_z:= \apmd (u\circ\psi^{-1})_{\psi(z)}\circ d\psi_z$$ is independent of the choice of chart. We say that $u$ is approximately metrically differentiable at $z$ and call $\apmd u_z$ the approximate metric derivative of $u$ at $z$. If $M$ is of dimension $n=1$ and $c\in W^{1,q}(M, X)$ then we abbreviate $$|c'(t)| = \apmd c_t(v),$$ where $v\in T_zM$ denotes either of the two unit vectors with respect to $g$.

Next, we specialize to the case that $M$ has dimension $n=2$. The $q$--energy of a semi-norm $s$ on (Euclidean) $\R^2$ is defined by  $$\mathbf{I}_+^q(s):= \max\{s(v)^q: v\in\R^2, |v|=1\}.$$
The jacobian of a norm $s$ on $\R^2$ is the unique number $\jac(s)$ such that $$\hm^2_{(\R^2, s)}(A) = \jac(s) \cdot |A|$$ for some and thus every subset $A\subset\R^2$ satisfying $|A|>0$. For a degenerate semi-norm $s$ we set $\jac(s):= 0$. If $s$ is a semi-norm on $\R^2$ and $L\colon\R^2\to\R^2$ is a linear map then 
\begin{equation}\label{eq:jac-semi-comp-linear}
\jac(s\circ L) = \jac(s)\cdot |\det L|.
\end{equation}
We define the jacobian and energy of a semi-norm $s$ on $(T_zM, g(z))$ by identifying it with $(\R^2,|\cdot|)$ via a linear isometry.  Notice that we always have $\jac(s)\leq \mathbf{I}_+^2(s)$.

\bd
 The Reshetnyak $q$-energy of $u\in W^{1,q}(\Omega, X)$ with respect to $g$ is defined by $$E_+^q(u, g):= \int_{\Omega} \mathbf{I}_+^q(\apmd u_z)\,d\hm^2_{g}(z).$$
\ed

We define the energy of the restriction of $u$ to a measurable subset $A\subset \Omega$ analogously. It can be shown that the $q$--th root of $E_+^q(u,g)$ is equal to the $L^q$--norm of the minimal weak upper gradient of $u$, see \cite[Theorem 7.1.20]{HKST15}.
If $q=2$ and $(U, \psi)$ is a conformal chart of $M$ then it follows from the area formula that $$E_+^2(u|_K, g) = \int_{\psi(K)} \mathbf{I}_+^2(\apmd(u\circ\psi^{-1})_w)\,dw = E_+^2(u\circ\psi^{-1}|_{\psi(K)}, g_{\rm Eucl})$$ for every compact set $K\subset U$, where $g_{\rm Eucl}$ denotes the Euclidean metric. This implies that the energy $E_+^2$ is invariant under precompositions with conformal diffeomorphisms.

\bd
 The parametrized (Hausdorff) area of $u\in W^{1,2}(\Omega, X)$ is defined by 
 \begin{equation}\label{eq:def-area}
  \Area(u):= \int_{\Omega} \jac(\apmd u_z)\,d\hm^2_g(z).
 \end{equation}
\ed

We notice that $\Area(u) \leq E_+^2(u,g)$ for all $g$. If $(U, \psi)$ is any chart of $M$ and $K\subset U$ is compact then $$\Area(u|_K) = \int_{\psi(K)} \jac(\apmd (u\circ\psi^{-1})_w)\,dw = \Area(u\circ\psi^{-1}|_{\psi(K)})$$ by \eqref{eq:jac-semi-comp-linear} and the area formula. As a consequence, the parametrized area of a Sobolev map is invariant under precompositions with biLipschitz homeo\-morphisms. Finally, if $u$ satisfies Lusin's property (N) then the area formula \cite{Kir94}, \cite{Kar07} for metric space valued maps yields $$\Area(u) = \int_X\#\{z\in \Omega: u(z) = x\}\,d\hm^2_X(x).$$

Next, we recall the definition of the trace of a Sobolev map. Let $\Omega \subset M\setminus \partial M$ be a Lipschitz domain. Then for every $z\in \partial \Omega$ there exist an open neighborhood $U\subset M$ and a biLip\-schitz map $\psi\colon (0,1)\times [0,1)\to M$ such that $\psi((0,1)\times (0,1)) = U\cap \Omega$ and $\psi((0,1)\times\{0\}) = U\cap \partial\Omega$. Let $u\in W^{1,q}(\Omega, X)$. Then for almost every $s\in (0,1)$ the map $t\mapsto u\circ\psi(s,t)$ has an absolutely continuous representative which we denote by the same expression. The trace of $u$ is defined by $$\trace(u)(\psi(s,0)):= \lim_{t\searrow 0} (u\circ\psi)(s,t)$$ for almost every $s\in(0,1)$. It can be shown (see \cite{KS93}) that the trace is independent of the choice of the map $\psi$ and defines an element of $L^q(\partial \Omega, X)$.

\bl\label{lem:bound-L2-norm}
 Let $M$ be a smooth compact and connected surface with non-empty boundary, and let $g$ be a Riemannian metric on $M$. Then there exists a constant $C$ with the following property. Let $(X, d)$ be a complete metric space, $x_0\in X$ and $R>0$. If $u\in W^{1,2}(M, X)$ satisfies $\trace(u)(z)\in \bar{B}(x_0, R)$ for almost every $z\in \partial M$ then $$\int_Md^2(x_0, u(z))\,d\hm_g^2(z)\leq C\cdot\left( R^2 + E_+^2(u, g)\right).$$
\el

\begin{proof}
 It suffices to prove the lemma in the case that $(M, g)$ has constant curvature $-1$, $0$, or $1$ and $\partial M$ is geodesic. Define a $1$--Lipschitz function $\psi$ on $X$ by $$\psi(x):= \max\{0, d(x, x_0) - R\}$$ and notice that the composition $v:= \psi\circ u$ belongs to $W^{1,2}(M\setminus\partial M)$ and satisfies $\trace(v)(z) = 0$ for almost every $z\in\partial M$. Since $$d^2(x_0, u(z)) \leq (v(z) + R)^2 \leq 2\cdot (v^2(z) + R^2)$$ for almost all $z\in M$ it is enough to show that there exists a constant $C$ depending on $(M,g)$ such that 
\begin{equation}\label{eq:Poincare-Sobolev-boundary}
 \int_M v^2\,d\hm_g^2\leq C\cdot E_+^2(u,g).
\end{equation}
In order to prove \eqref{eq:Poincare-Sobolev-boundary}, let $(\hat{M}, \hat{g})$ be the Schottky double of $(M, g)$ obtained by gluing two copies $M_-$ and $M_+$ of $M$ along their boundaries and by doubling $g$, see for example \cite[Chapter 4.4]{Jos06-cptRiem} or \cite[Chapter II.1.3]{AS60}. Since $(M, g)$ has constant curvature $-1$, $0$, or $1$ and $\partial M$ is geodesic the Schottky double $(\hat{M}, \hat{g})$ is a smooth compact surface without boundary of the same curvature.
Define a map $\hat{v}\colon\hat{M}\to\R$ by $$\hat{v}(z):= \left\{\begin{array}{rl}
 v(z) & z\in M_+\\
 -v(z) & z\in M_-.
\end{array}\right.
$$
Since $v$ has zero trace it follows that $\hat{v}$ belongs to $W^{1,2}(\hat{M})$. Because $$\int_{\hat{M}} \hat{v}\,d\hm_{\hat{g}}= 0$$ 
the Sobolev-Poincar\'e inequality \cite[Theorem 2.11]{Heb99} implies that
\begin{equation*}
  \int_M v^2\,d\hm_g^2= \frac{1}{2}\cdot \int_{\hat{M}}|\hat{v}|^2\,d\hm_{\hat{g}}^2\leq C'\cdot\left(\int_{\hat{M}} |\nabla \hat{v}|\,d\hm^2_{\hat{g}}\right)^2= 4C'\cdot\left(\int_M|\nabla v|\,d\hm_g^2\right)^2
\end{equation*}
for some constant $C'$ depending on $(M,g)$. H\"older's inequality yields 
\begin{equation*}
  \left(\int_M|\nabla v|\,d\hm_g^2\right)^2 \leq \hm_g^2(M)^2\cdot \int_{M}|\nabla v|^2\,d\hm_g^2\leq \hm_g^2(M)^2\cdot E_+^2(u,g),
 \end{equation*}
  which, together with the inequality above, establishes \eqref{eq:Poincare-Sobolev-boundary}. This completes the proof of the lemma.
\end{proof}

\section{Relative systole, hyperbolic collars, and Mumford compactness}\label{sec:rel-sys}

We denote by $\Hyp$ the hyperbolic plane. For us, it will be most convenient to work with the upper-half plane model of hyperbolic space, so we set $$\Hyp:= \{z=x+iy\in \C: y>0\}$$ and equip $\Hyp$ with the Riemannian metric $g_{\Hyp}:= \frac{1}{y^2}\cdot (dx^2 + dy^2)$.

Let $M$ be a smooth compact and connected surface with non-empty boundary and of strictly negative Euler characteristic. Then there exists a Riemannian metric $g$ on $M$ such that $(M,g)$ has constant curvature $-1$ and $\partial M$ is geodesic. Any such metric will be called hyperbolic metric on $M$. 

We will use the following variant of the systole adapted to the setting of surfaces with boundary. For the next definition and the proposition below we fix a hyperbolic metric $g$ on $M$.

\bd
The relative systole $\sys_{\rm rel}(M,g)$ of $(M,g)$ is the minimal length of curves $\gamma$ in $M$ of the following form. Either $\gamma$ is closed and is not contractible in $M$ via a family of closed curves. Or the endpoints of $\gamma$ lie on the boundary $\partial M$ of $M$ and $\gamma$ is not contractible via a family of curves with endpoints on $\partial M$. 
\ed

Using the relative systole we now state a simple consequence of the well-known collar lemma from hyperbolic geometry. 
Let $\lambda_0>0$ satisfy 
\begin{equation}\label{eq:lambda-0}
 \sinh(\lambda_0) \cdot \sinh(2) = 1
\end{equation}
and suppose the relative systole $\lambda:= \sys_{\rm rel}(M,g)$ satisfies $\lambda\leq \lambda_0$. Let $\gamma$ be a curve as in the definition above of minimal length. Thus, $\gamma$ is a geodesic of length $\lambda$, and we parametrize it by arc-length on the interval $[0,\lambda]$. Let $A_\lambda$ be the subset of $\Hyp$ given by 
\begin{equation}\label{eq:A-lambda}
 A_\lambda:= \left\{ e^{s+it}: 0\leq s\leq \lambda, \frac{\pi}{4}\leq t\leq \frac{\pi}{2}\right\}
\end{equation}
 and let $\Sigma_\lambda$ be the hyperbolic surface obtained from $A_\lambda$ by identifying $e^{it}$ with $e^{\lambda + it}$ for all $t$.

\bp\label{prop:map-alpha}
 There exists a smooth map $\alpha\colon A_\lambda\to (M,g)$ with the following properties:
 \begin{enumerate}
  \item[(i)]  For all $s\in[0,\lambda]$ we have $\alpha(e^{s+i\frac{\pi}{2}}) = \gamma(s)$.
  \item[(ii)] If $\gamma$ is a closed curve then $\alpha$ descends to an isometric map $\Sigma_\lambda\to (M,g)$. 
  \item[(ii)] If $\gamma$ has distinct endpoints then $\alpha$ is isometric and maps the geodesic segments $\{e^{it}: \frac{\pi}{4}\leq t\leq \frac{\pi}{2}\}$ and $\{e^{\lambda+ it}: \frac{\pi}{4}\leq t\leq \frac{\pi}{2}\}$ to geodesic segments on $\partial M$.
 \end{enumerate}
\ep

\begin{proof}
Consider the Schottky double $(\hat{M}, \hat{g})$ of $(M, g)$, see the proof of Lemma~\ref{lem:bound-L2-norm}, which is a closed hyperbolic surface. If $\gamma$ is a closed curve then it defines a simple closed geodesic in $(\hat{M}, \hat{g})$ of length $\lambda$. If $\gamma$ is not closed then it intersects $\partial M$ perpendicularly and thus the doubled curve defines a simple closed geodesic in $(\hat{M}, \hat{g})$ of length $2\lambda$, see \cite{Bus10}.

The existence of a map $\alpha$ as in the proposition now follows from the collar lemma (see \cite[Theorem 4.1.1]{Bus10}), taking into account that the width of $A_\lambda$ (i.e. the length $\xi$ of the geodesic segment $\{e^{it}: \frac{\pi}{4}\leq t\leq \frac{\pi}{2}\}\subset A_\lambda$) is bounded by $$\xi = \int_{\frac{\pi}{4}}^{\frac{\pi}{2}} \frac{1}{\sin t}\,dt \leq \int_{\frac{\pi}{4}}^{\frac{\pi}{2}} \frac{1}{\sin^2t}\,dt=1,$$ as well as that $\lambda\leq \lambda_0$ and $\lambda_0$ satisfies \eqref{eq:lambda-0}.
\end{proof}

We moreover need a variant of the Mumford compactness theorem for surfaces with non-empty boundary. 

\bt\label{thm:Mumford-cptness-boundary}
 Let $(g_n)$ be a sequence of hyperbolic metrics on $M$. If $$\inf\{ \sys_{\rm rel}(M, g_n): n\in\N\}>0$$ then there exist diffeomorphisms $\varphi_n\colon M \to M$ such that a subsequence of $(\varphi_n^* g_n)$ converges smoothly to a hyperbolic metric on $M$. 
\et

\begin{proof}
Let $(\hat{M}_n, \hat{g}_n)$ be the Schottky double of $(M, g_n)$. Notice that $\hat{g}_n$ is invariant under the natural involution on $\hat{M}_n$. We identify $\hat{M}_n$ with $\hat{M}:= \hat{M}_1$ via a diffeomorphism which commutes with the involution and pull back the metric $\hat{g}_n$ to $\hat{M}$ via this diffeomorphism. 

Since $\sys_{\rm rel}(M, g_n)$ bounds from below the usual systole of $(\hat{M}, \hat{g}_n)$ it follows from the Mumford compactness theorem (see \cite[Theorem 4.4.1]{DHT10}) that there exist diffeomorphisms $\hat{\varphi}_n\colon \hat{M} \to \hat{M}$ which commute with the involution $\iota\colon \hat{M}\to \hat{M}$, map each half of $\hat{M}$ to itself and such that, after possibly passing to a subsequence, the Riemannian metrics $\hat{\varphi}_n^* \hat{g}_n$ converge smoothly to a hyperbolic metric $\hat{g}$ on $\hat{M}$. Since $\hat{g}$ is invariant under $\iota$ it follows that the fixed point set $\partial M$ of $\iota$ is totally geodesic with respect to $\hat{g}$. The theorem now follows from restricting to one of the two halves of $\hat{M}$.
\end{proof}

\section{Energy minimizers are infinitesimally isotropic}\label{sec:energy-min-isotropic}

Let $M$ be a smooth compact surface and $X$ a complete metric space. The following definition appears in \cite{LW17-en-area} and implicitly in \cite{LW15-Plateau} in the case that $(M,g)$ is a two-dimensional Euclidean domain.

\bd
 A map $u\in W^{1,2}(M, X)$ is infinitesimally isotropic with respect to a Riemannian metric $g$ if for almost every $z\in M$ the approximate metric derivative $\apmd u_z$ is either zero or it is a norm and the ellipse of maximal area contained in the unit ball with respect to $\apmd u_z$ is a round ball with respect to $g$.
\ed

Suppose $u$ is infinitesimally isotropic with respect to $g$. Then for almost every $z\in M$ the approximate metric derivative $s:= \apmd u_z$ satisfies $\mathbf{I}_+^2(s) \leq \frac{4}{\pi}\cdot \jac(s)$ by \cite[Theorem 6.2]{Bal97}. In particular it follows that 
\begin{equation}\label{eq:energy-bdd-by-area-inf-isotropic}
E_+^2(u,g) \leq \frac{4}{\pi}\cdot \Area(u).
\end{equation}
 Moreover, $u$ is $Q$--quasiconformal with respect to $g$  with $Q=\sqrt{2}$ in the sense that for almost all $z\in M$ we have $$\apmd u_z(v) \leq Q\cdot \apmd u_z(w)$$ for all $v, w\in T_zM$ with $|v|_z=|w|_z$. Here, $|v|_z$ is the length of $v$ with respect to $g$. If $X$ is a Riemannian manifold or, more generally, if $X$ has property (ET) in the sense of \cite{LW15-Plateau} and $u$ is infinitesimally isotropic then $u$ is $Q$--quasiconformal with $Q=1$. 

The following result generalizes \cite[Theorem 6.2]{LW15-Plateau} from the setting of Sobolev maps defined on the Euclidean unit disc to that of Sobolev maps defined on a smooth compact surface. 

\bt\label{thm:energy-min-isotropic}
 Let $M$ be a smooth compact and connected surface equipped with a Riemannian metric $g$. Let $X$ be a complete metric space and $u\in W^{1,2}(M, X)$. If for every Riemannian metric $g'$ on $M$ and every biLipschitz homeomorphism $\varphi\colon (M, g')\to (M, g)$ we have $$E_+^2(u, g)\leq E_+^2(u\circ\varphi, g')$$ then $u$ is infinitesimally isotropic with respect to $g$.
\et

As the proof will show, it is enough to consider Riemannian metrics $g'$ of constant curvature $-1$, $0$, or $1$ such that $\partial M$ is geodesic.

\begin{proof}
 We will show that every point in $M\setminus \partial M$ has an open neighborhood on which $u$ is infinitesimally isotropic with respect to $g$. 
 
 Fix $x\in M\setminus \partial M$ and let $(U_0,\psi_0)$ be a smooth chart around $x$ with the properties that $\partial U_0$ is smooth and $U_0\cap \partial M=\emptyset$, the chart map $\psi_0$ is conformal, orientation preserving and biLipschitz and satisfies $\psi_0(U_0) = D$. We will show that for almost every $z\in D$ we have 
 \begin{equation}\label{eq:I+composed-SL}
 \mathbf{I}_+^2(\apmd(u\circ\psi_0^{-1})_z) \leq \mathbf{I}_+^2(\apmd(u\circ\psi_0^{-1})_z\circ T)
 \end{equation} 
 for every $T\in\SL_2(\R)$. From the proof of \cite[Lemma 6.5]{LW15-Plateau} it then follows that $u\circ\psi_0^{-1}$ is infinitesimally isotropic on $D$ with respect to the Euclidean metric $g_{\rm Eucl}$ and thus $u$ is infinitesimally isotropic on $U_0$ with respect to $g$ by the conformality of $\psi_0$. 
 
 We argue by contradiction and assume that \eqref{eq:I+composed-SL} fails on a set of strictly positive measure. The proof of \cite[Theorem 6.2]{LW15-Plateau} then shows that there exists a biLipschitz homeomorphism $\rho\colon D\to D$ which satisfies 
 \begin{equation}\label{eq:energy-lower-disc}
 E_+^2(u\circ\psi_0^{-1}\circ\rho^{-1}, g_{\rm Eucl}) < E_+^2(u\circ\psi_0^{-1}, g_{\rm Eucl})
 \end{equation}
  and which is smooth and conformal outside some compact subset $B\subset D$ with smooth boundary. We may of course assume that $\rho$ is orientation preserving.
 
Next, we complete the chart $(U_0, \psi_0)$ to a finite atlas $\mathcal{A} = \{(U_i, \psi_i): i=0,\dots, n\}$ of $M$ consisting of conformal charts which are all orientation preserving and such that for any $i\not= 0$ the set $U_i$ does not intersect $B':=\psi_0^{-1}(B)$. On the topological manifold $M$ we consider a new atlas $\mathcal{A}'$ given by $$\mathcal{A}' = \{ (U'_i, \psi'_i): i=0,\dots, n\},$$ where $U'_i = U_i$ for all $i$ and $\psi'_0 = \rho\circ \psi_0$ and $\psi'_i=\psi_i$ for $i\not=0$. Since the transition maps $\psi'_i\circ (\psi'_j)^{-1}$ are conformal the atlas $\mathcal{A}'$ induces a new conformal (and thus also a new smooth) structure on $M$. We denote by $M'$ the resulting Riemann surface. By the uniformization theorem there exists a Riemannian metric $g'$ on $M'$ of constant curvature $-1$, $0$, or $1$ for which all the charts $(U'_i, \psi'_i)$ are conformal and such that the boundary $\partial M'$ is geodesic.

Let $\varphi\colon (M',g')\to (M,g)$ be the identity map. Since $$\psi_j\circ\varphi\circ(\psi'_0)^{-1} = (\psi_j\circ\psi_0^{-1})\circ \rho^{-1}$$ and $\psi_j\circ\varphi\circ (\psi'_i)^{-1} = \psi_j\circ \psi_i^{-1}$ when $i\not=0$ it follows from the properties of $\psi_j$ and $\rho$ that $\varphi$ is smooth and conformal on $M'\setminus B'$ and locally biLipschitz on $U'_0$. From this we see that $\varphi$ is a biLipschitz homeomorphism. Moreover, since $\partial U'_0$ is smooth in $M'$ and hence has zero measure it follows from \eqref{eq:energy-lower-disc} and the conformal invariance of the energy that
\begin{equation*}
 \begin{split}
  E_+^2(u\circ\varphi, g') &=  E_+^2(u\circ\varphi|_{U'_0}, g') + E_+^2(u\circ\varphi|_{M'\setminus \overline{U'_0}}, g')\\
  &= E_+^2(u\circ(\rho\circ\psi_0)^{-1}, g_{\rm Eucl}) + E_+^2(u|_{M\setminus U_0}, g)\\
  &<E_+^2(u\circ\psi_0^{-1}, g_{\rm Eucl}) + E_+^2(u|_{M\setminus U_0}, g)\\
  &= E_+^2(u, g).
  \end{split}
 \end{equation*}
Since $M'$ and $M$ are homeomorphic and thus also diffeomorphic we may assume, after pulling back $g'$ by a diffeomorphism, that $M'=M$ and thus the above contradicts our hypotheses. This shows that $u$ is infinitesimally isotropic on $U_0$ with respect to $g$ and concludes the proof.
\end{proof}

\section{Equi-continuity of traces}\label{sec:equi-cont-traces}

Let $(X,d)$ be a complete metric space and let $\Gamma\subset X$ be the disjoint union of $k\geq 1$ rectifiable Jordan curves $\Gamma_1,\dots, \Gamma_k$. Let furthermore $M$ be a smooth compact and connected surface with $k$ boundary components and of genus $p\geq 0$ such that $k+2p\geq 2$. Thus, $M$ is not diffeomorphic to $\overline{D}$. We assume furthermore that $X$ admits a local quadratic isoperimetric inequality. 

\bp\label{prop:equi-cont-traces}
 Let $g$ be a Riemannian metric on $M$. Then for every $\eta>0$ and $K>0$ the family $$\{\trace(u): u\in\Lambda(M, \Gamma, X), E_+^2(u, g)\leq K, \Area(u) \leq a^*_p(\Gamma, X) - \eta\}$$ is equi-continuous.
\ep

In the above, $\trace(u)$ refers to the continuous representative of the trace of $u$. When $M$ is diffeomorphic to $\overline{D}$ an analogous statement holds. In this case the condition on the area in the statement is replaced by a $3$--point condition, see \cite{LW15-Plateau}.

In the sequel we will use the following terminology. A smooth compact surface $M^*$ with $k$ boundary components is called a reduction of $M$ if $M^*$ has the following property. Either $M^*$ is connected and has genus at most $p-1$ or $M^*$ has exactly two connected components, each of which has non-empty boundary, and the genus of $M^*$ is $p$. We thus have $$a_p^*(\Gamma, X) = \min\{a(M^*, \Gamma, X): \text{ $M^*$ is a reduction of $M$}\}$$ by the definition of $a_p^*(\Gamma, X)$.

\begin{proof}[Proof of Proposition~\ref{prop:equi-cont-traces}]
 Let $C, l_0>0$ be the constants from the local quadratic isoperimetric inequality, see Definition~\ref{def:loc-quad-isop}.
 Let $\eta, K>0$ and $\varepsilon>0$ and set $$\rho:= \min\left\{\varepsilon, \frac{l_0}{2}, \sqrt{\frac{\eta}{8C}}\right\}.$$
 There exists $0<\rho'<\rho$ such that whenever $x,x'\in\Gamma$ are distinct points belonging to the same Jordan curve $\Gamma_j$ and satisfying $d(x,x')<\rho'$ then the shorter of the two subcurves of $\Gamma_j$ connecting $x$ and $x'$ has length at most $\rho$.
Let $0<\delta<1$ be so small that $$\pi\cdot \left(\frac{8K}{|\log(\delta)|}\right)^{\frac{1}{2}} < \rho'$$ and that every point $z_0\in \partial M$ has a neighborhood in $M$ which is the image of the set $$B:=\left\{z\in\C: \text{$|z|\leq 1$ and $|z-1|<\sqrt{\delta}$}\right\}$$ under a diffeomorphism $\psi$ that is $2$--biLipschitz and maps the point $1\in B$ to $z_0$.

Let $u$ be a map from the family defined in the statement of the proposition and let $\gamma$ be a component of $\partial M$. It suffices to show that the continuous representative of $\trace(u)$ maps segments in $\gamma$ of length at most $2\delta$ to curves of length at most $\varepsilon$. Fix $z_0\in\partial M$ and let $\psi$ be a diffeomorphism as above. Let $\eta_r$ be a constant speed parametrization (defined on some interval $I$) of the curve $\{z\in B: |z-1|= r\}$ whenever $r\in(0,\sqrt{\delta})$ and set $\beta_r:= \psi\circ \eta_r$. 

By the Courant-Lebesgue lemma (see \cite[Lemma 7.3]{LW15-Plateau}) there exists a set $A\subset (\delta, \sqrt{\delta})$ of strictly positive measure such that for every $r\in A$ the map $u\circ\beta_r$ is in $W^{1,2}(I,X)$ and its absolutely continuous representative, again denoted by $u\circ\beta_r$, satisfies $$\length(u\circ\beta_r)  \leq \pi\cdot \left(\frac{2 E_+^2(u\circ \psi, g_{\rm Eucl})}{|\log(\delta)|}\right)^{\frac{1}{2}} \leq \pi\cdot \left(\frac{8 E_+^2(u, g)}{|\log(\delta)|}\right)^{\frac{1}{2}}< \rho'.$$ 
For almost every $r\in A$ the endpoints of $u\circ\beta_r$ coincide with $\trace(u)(a_r)$ and $\trace(u)(b_r)$, where $a_r$ and $b_r$ are the endpoints of $\beta_r$, and in particular
\begin{equation}\label{eq:endpoints-close-equicont}
 d(\trace(u)(a_r), \trace(u)(b_r))\leq \length(u\circ\beta_r)<\rho'
\end{equation}
 by the above. Moreover, for almost every $r\in A$ we have $$\trace(u\circ\psi|_{B_r})\circ\eta_r=u\circ\beta_r,$$ where $B_r:= \{z\in \C: \text{$|z|<1$ and $|z-1|<r$}\}$. Fix $r\in A$ such that all of the above hold.

Let $\gamma^-$ be the subcurve of $\gamma$ which connects $a_r$ and $b_r$ and contains $z_0$. Let $j$ be such that $\trace(u)(\gamma) = \Gamma_j$ and denote by $\Gamma_j^-$ the image of $\gamma^-$ under $\trace(u)$. If we can show that $\Gamma_j^-$ has length at most $\rho$ (and thus no larger than $\varepsilon$) then the proof of equi-continuity is complete.
We argue by contradiction and assume that $\length(\Gamma_j^-)>\rho$. 

Let $U\subset M$ be the Jordan domain enclosed by the concatenation $\beta_r \cup \gamma^-$, thus $U = \psi(B_r)$ and $\gamma^-$ corresponds to $\psi(\{z\in \C: \text{$|z|= 1$ and $|z-1|\leq r$}\})$. 
Let $\gamma^+$ be the complementary segment in $\gamma$ and denote its image under $\trace(u)$ by $\Gamma_j^+$. Since $\length(\Gamma_j^-)>\rho$ it follows from \eqref{eq:endpoints-close-equicont} and the choice of $\rho'$ that $\length(\Gamma_j^+)\leq \rho$. In particular, the closed curve $(u\circ\beta_r)\cup \Gamma_j^+$ has length bounded by $$\length(u\circ\beta_r) + \length(\Gamma_j^+)< \rho' + \rho<2\rho.$$ 
Set $\Omega^-:= D\setminus\overline{B_r}$. Since $u\circ\beta_r$ is a $W^{1,2}$--curve it follows from \cite[Lemma 8.5]{LW15-Plateau} that there exists $w^-\in W^{1,2}(\Omega^-, X)$ with $$\Area(w^-)< 4C\cdot \rho^2,$$ such that $\trace(w^-)\circ\eta_r = u\circ\beta_r$ and the restriction of $\trace(w^-)$ to $S^1\setminus \overline{B_r}$ is a constant speed parametrization of $\Gamma_j^+$. By the gluing theorem \cite[Theorem 1.12.3]{KS93} the map $w$ coinciding with $u\circ\psi$ on $B_r$ and with $w^-$ on $\Omega^-$ defines an element of $W^{1,2}(D, X)$. Clearly, we have $$\Area(w) = \Area(u\circ\psi|_{B_r}) + \Area(w^-)< \Area(u|_U) + 4C\rho^2$$ and the trace of $w$ is a weakly monotone parametrization of $\Gamma_j$. 

Suppose first that $k=1$. Then $M^*:= \overline{D}$ is a reduction of $M$ and $w\in \Lambda(M^*, \Gamma, X)$. Moreover, we have $$a_p^*(\Gamma, X)\leq \Area(w) < \Area(u) + \eta,$$ which is impossible.

Suppose now that $k\geq 2$. Let $M'$ be the smooth surface obtained from $M$ by gluing a disc to the boundary component $\gamma$ of $\partial M$ and view $M$ as a subset of $M'$. Let $\Omega^+\subset M'$ be the Jordan domain enclosed by $\gamma^+\cup \beta_r$. Then $\gamma^+\cup \beta_r$ is the common boundary of the two disjoint Lipschitz domains $\Omega^+$ and $M \setminus (\overline{U}\cup \gamma)$. Since the continuous representative of $\trace(u|_{M\setminus \overline{U}})|_{\gamma^+\cup\beta_r}$ has length less than $2\rho$ and is the trace of a Sobolev annulus it follows from the local quadratic isoperimetric inequality and from \cite[Lemma 4.8]{LW-intrinsic} and its proof that there exists $w^+\in W^{1,2}(\Omega^+, X)$ whose trace coincides with $\trace(u|_{M\setminus \overline{U}})$ on the common boundary $\partial\Omega^+$ and whose area satisfies $$\Area(w^+)< 4C\rho^2.$$
The gluing theorem \cite[Theorem 1.12.3]{KS93} shows that the map $u'$ agreeing with $u$ on $M\setminus \overline{U}$ and with $w^+$ on $\Omega^+$ defines an element of $W^{1,2}(M', X)$ and $$\Area(u')= \Area(u|_{M\setminus \overline{U}}) + \Area(w^+)< \Area(u|_{M\setminus \overline{U}}) + 4C\rho^2.$$

Finally, let $M^*$ be the disjoint union of $\overline{D}$ and $M'$. Then $M^*$ has $k$ boundary curves, has genus $p$ and has two connected components each of which has non-empty boundary. In particular, $M^*$ is a reduction of $M$. The map $v\in W^{1,2}(M^*, X)$ which agrees with $w$ on $D$ and with $u'$ on $M'$ belongs to $\Lambda(M^*, \Gamma, X)$ and satisfies $$\Area(v) = \Area(w) + \Area(u') < \Area(u) + 8C\rho^2\leq \Area(u) + \eta.$$ This is impossible since $\Area(u)\leq a^*_p(\Gamma, X) - \eta\leq \Area(v) - \eta$. We conclude that the length of $\Gamma_j^-$ is no larger than $\rho$ and thus no larger than $\varepsilon$, finishing the proof of the proposition.
\end{proof}

\section{A lower bound for the relative systole}\label{sec:lower-bound-rel-sys}

Let $(X,d)$ be a complete metric space and let $\Gamma\subset X$ be the disjoint union of $k\geq 1$ rectifiable Jordan curves. Let furthermore $M$ be a smooth compact and connected surface with $k$ boundary components and of genus $p\geq 0$. We assume that $X$ admits a local quadratic isoperimetric inequality. We furthermore assume that $k+2p\geq 3$, that is, $M$ has strictly negative Euler characteristic and hence admits a hyperbolic metric.

\bp\label{prop:bound-systole}
For every $\eta>0$ and $K>0$ there exists $\varepsilon>0$ with the following property. Let $u\in \Lambda(M, \Gamma, X)$ be such that $$\Area(u) \leq a^*_p(\Gamma, X) - \eta$$ and let $g$ be a hyperbolic metric on $M$.  If $E_+^2(u, g)\leq K$ then $\sys_{\rm rel}(M,g)$ is bounded from below by $\varepsilon$.
\ep

The proof of Proposition~\ref{prop:bound-systole} relies on the following lemma. Let $\lambda_0>0$ satisfy \eqref{eq:lambda-0} and let $g$ be a hyperbolic metric on $M$ with $$\lambda:=\sys_{\rm rel}(M, g)\leq \lambda_0.$$ Let $\alpha\colon A_\lambda\to M$ be the map from Proposition~\ref{prop:map-alpha} for some curve $\gamma$ in $M$ realizing the relative systole of $(M,g)$, where the set $A_\lambda\subset\Hyp$ is as in \eqref{eq:A-lambda}. Set $I:= [0,\lambda]$ and let $\alpha_t\colon I\to M$ be the curve defined by $\alpha_t(s):= \alpha(e^{s+it})$ whenever $\frac{\pi}{4}<t< \frac{\pi}{2}$. If $\gamma$ is a closed curve then $\alpha_t$ defines a smooth closed curve for every $t$. In this case we interpret $\alpha_t$ as a map from the smooth one-dimensional manifold obtained from $I$ by identifying its endpoints. We will denote this again by $I$.

\bl\label{lem:short-image-curves-collar}
 Given $u\in W^{1,2}(M, X)$ there exists a subset $A\subset (\frac{\pi}{4}, \frac{\pi}{2})$ of strictly positive measure such that $u\circ\alpha_t\in W^{1,2}(I, X)$ for every $t\in A$, and the absolutely continuous representative of $u\circ\alpha_t$ satisfies 
 \begin{equation}\label{eq:length-u-alpha-t} 
  \length(u\circ \alpha_t) \leq 2\left[\lambda \cdot E_+^2(u,g)\right]^{\frac{1}{2}}.
 \end{equation}
\el

\begin{proof}
 Define a biLipschitz map $\varphi\colon I\times [\frac{\pi}{4}, \frac{\pi}{2}]\to A_\lambda$ by $\varphi(s,t):= e^{s+it}$ and notice that $\alpha_t(s) = (\alpha\circ\varphi)(s,t)$. 
 By \cite[Proposition 4.10]{LW15-Plateau} and its proof we have that $u\circ\alpha_t\in W^{1,2}(I, X)$ and $$|(u\circ\alpha_t)'(s)| = \apmd (u\circ\alpha)_{\varphi(s,t)} \left(\frac{\partial \varphi}{\partial s}(s,t)\right)$$ for almost every $t$ and $s$. In particular, 
 $$|(u\circ\alpha_t)'(s)| \leq \mathbf{I}_+^1(\apmd (u\circ\alpha)_{\varphi(s,t)})\cdot \frac{1}{\sin t},$$ and H\"older's inequality implies that the absolutely continuous representative of $u\circ\alpha_t$, denoted by the same expression, satisfies $$\length^2(u\circ\alpha_t)= \left(\int_0^{\lambda} |(u\circ\alpha_t)'(s)|\, ds\right)^2 \leq \frac{\lambda}{\sin^2 t}\cdot \int_0^{\lambda} \mathbf{I}_+^2(\apmd (u\circ\alpha)_{\varphi(s,t)})\,ds$$ for almost every $t$. Since $\alpha$ is locally isometric and injective in the interior of $A_\lambda$ and since $|\det(d\varphi_{(s,t)})| = \frac{1}{\sin^2 t}$ we obtain, using the area formula, that $$\int_{\frac{\pi}{4}}^{\frac{\pi}{2}} \length^2(u\circ\alpha_t)\, dt \leq \lambda\cdot E_+^2(u\circ\alpha, g_{\Hyp})\leq \lambda\cdot E_+^2(u, g).$$ Thus, inequality \eqref{eq:length-u-alpha-t} cannot fail for almost every $t$, completing the proof.
\end{proof}

In the following proof we use the notation introduced at the beginning of this section and before the statement of Lemma~\ref{lem:short-image-curves-collar}.

\begin{proof}[Proof of Proposition~\ref{prop:bound-systole}]
 Let $C, l_0>0$ be the constants from the local quadratic isoperimetric inequality, see Definition~\ref{def:loc-quad-isop}.
 Let $\eta>0$ and $K>0$ and set $$\rho:= \min\left\{\frac{l_0}{2}, \sqrt{\frac{\eta}{8C}}\right\}.$$ Recall that $\Gamma$ is the disjoint union of finitely many rectifiable Jordan curves $\Gamma_j$. Therefore, there exists $0<\rho'<\rho$ such that whenever $x,x'\in\Gamma$ are distinct points with $d(x,x')<\rho'$ then they belong to a single Jordan curve $\Gamma_j$ and the shorter of the two subcurves of $\Gamma_j$ connecting $x$ and $x'$ has length at most $\rho$. We will show that the proposition holds with $$\varepsilon:= \min\left\{\frac{\rho'^2}{4K}, \lambda_0\right\}.$$
Let $g$ be a hyperbolic metric on $M$. Let $u\in \Lambda(M, \Gamma, X)$ and suppose that $$\Area(u)\leq a^*_p(\Gamma, X) - \eta$$ and $E_+^2(u, g)\leq K$. 
We argue by contradiction and assume that $\lambda:=  \sys_{\rm rel}(M, g)$ satisfies $\lambda<\varepsilon$. By Lemma~\ref{lem:short-image-curves-collar} and the choice of $\varepsilon$ there exists a subset $A\subset (\frac{\pi}{4}, \frac{\pi}{2})$ of strictly positive measure such that $u\circ\alpha_t\in W^{1,2}(I, X)$ for all $t\in A$ and the absolutely continuous representative of $u\circ\alpha_t$ satisfies
 \begin{equation}\label{eq:length-short-parallel}
  \length(u\circ \alpha_t) \leq 2\sqrt{\lambda K} < 2\sqrt{\varepsilon K} \leq \rho'.
 \end{equation} 
 
 We distinguish two cases and first assume that $\alpha_t$ is a smooth closed curve for every $t$. Let $t\in A$ and let $M'$ be the smooth surface obtained by cutting $M$ along $\alpha_t$, so that $\alpha_t$ gives rise to two new boundary components in $M'$, denoted by $\alpha_t^-$ and $\alpha_t^+$. Embed $M'$ diffeomorphically into a smooth compact surface $M^*$ such that $M^*\setminus M'$ is the disjoint union $\Omega^-\cup\Omega^+$ of two open discs, with $\Omega^\pm$ being bounded by $\alpha_t^\pm$. Then $M'$ is (the closure of) a Lipschitz domain in $M^*$. Consider $u$ as an element of $W^{1,2}(M', X)$. We may assume that $\trace(u)\circ\alpha_t^\pm$ coincides with a copy of $u\circ\alpha_t$ (this is true for almost every $t\in A$). Since $u\circ\alpha_t$ is a $W^{1,2}$--curve whose length satisfies \eqref{eq:length-short-parallel}, it follows from \cite[Lemma 8.5]{LW15-Plateau} that there exist maps $w^\pm\in W^{1,2}(\Omega^\pm, X)$ whose traces coincide with $\trace(u)\circ\alpha_t^\pm$ and which satisfy $$\Area(w^\pm) \leq C\cdot \length^2(u\circ\alpha_t)< C\cdot \rho'^2<\frac{\eta}{2}.$$
By the gluing theorem \cite[Theorem 1.12.3]{KS93} the map $v$ coinciding with $u$ on $M'$ and with $w^\pm$ on $\Omega^\pm$ defines a Sobolev map on all of $M^*$ and the trace of $v$ agrees with that of $u$. Thus we have $v\in \Lambda(M^*, \Gamma, X)$ and moreover 
\begin{equation}\label{eq:reduction-smaller-area}
\Area(v)= \Area(u) + \Area(w^-)+\Area(w^+) < \Area(u) + \eta.
\end{equation}
 Finally, notice that the Euler characteristic $\chi(M^*)$ of $M^*$ satisfies $\chi(M^*) = \chi(M) + 2$. Therefore, if $M^*$ is connected then it has genus $p-1$ and is thus a reduction of $M$ in the sense defined after Proposition~\ref{prop:equi-cont-traces}. In this case we obtain with \eqref{eq:reduction-smaller-area} that $$a_p^*(\Gamma, X)\leq \Area(v)<\Area(u)+\eta,$$ which contradicts the assumption on $u$. Now, suppose that $M^*$ has two connected components. Then $M^*$ has genus $p$. If both connected components have non-empty boundary then $M^*$ is again a reduction and we obtain a contradiction exactly as above. If one of the components of $M^*$ has empty boundary then this component has genus at least one. Then the surface obtained by omitting this component is connected and has genus at most $p-1$; therefore it is a reduction of $M$ and we obtain a contradiction as before by considering the restriction of $v$ to this surface.

Let us turn to the second case and therefore assume that $\alpha_t$ has endpoints on $\partial M$ for all $t$, see Proposition~\ref{prop:map-alpha}. Denote by $a_t, b_t\in \partial M$ the endpoints of $\alpha_t$ and notice that for almost every $t\in A$ the endpoints of $u\circ\alpha_t$ coincide with $\trace(u)(a_t)$ and $\trace(u)(b_t)$. For such $t$ the inequality \eqref{eq:length-short-parallel} implies that $$d(\trace(u)(a_t), \trace(u)(b_t)) \leq \length(u\circ\alpha_t) < \rho'.$$ By the choice of $\rho'$,  the points $\trace(u)(a_t)$ and $ \trace(u)(b_t)$ thus lie on a single Jordan curve $\Gamma_j$ and the shorter subcurve of $\Gamma_j$ connecting the two points has length at most $\rho$. In particular, the points $a_t$ and $b_t$ lie on the same component of $\partial M$. Let $\gamma^-$ and $\gamma^+$ be the two segments in $\partial M$ connecting $a_t$ and $b_t$ and let $\Gamma_j^-$ and $\Gamma_j^+$ be their images under $\trace(u)$. If $\Gamma_j^-$ denotes the shorter of the two subcurves then $\length(\Gamma_j^-)\leq\rho$ by the above.

Let $M'$ be the manifold with corners obtained by cutting $M$ along the curve $\alpha_t$. Thus, $\alpha_t$ gives rise to two boundary curves in $M'$, which we denote by $\alpha_t^-$ and $\alpha_t^+$, so that $\gamma^-\cup \alpha_t^-$ and $\gamma^+\cup\alpha_t^+$ are piecewise smooth closed curves. Next, embed $M'$ diffeomorphically into a smooth compact surface $M^*$ with boundary such that $M^*\setminus M'$ is the disjoint union $\Omega^-\cup\Omega^+$ of open sets, where $\Omega^-$ is an open disc bounded by $\gamma^-\cup \alpha_t^-$ and $\Omega^+$ is an open half-disc bounded by the curve $\alpha_t^+$. Thus, $\Omega^+$ is diffeomorphic to $$D^+:= \{z\in \C: \text{$|z|<1$ and $\operatorname{Im}(z)\geq 0$}\}$$ and $\alpha_t^+$ corresponds to the subset $D^+\cap \{|z|= 1\}$. Let $\tilde{\gamma}^+$ be the subset of $\Omega^+$ corresponding to $D^+\cap \{\operatorname{Im}(z)=0\}$. Then $\tilde{\gamma}^+$ and the curve $\gamma^+$ together form one component of $\partial M^*$.

 Consider $u$ as an element of $W^{1,2}(M', X)$. Exactly as above, we may assume that $\trace(u)\circ\alpha_t^\pm$ coincides with a copy of $u\circ\alpha_t$. Since $$\length(\Gamma_i^-) + \length(u\circ\alpha_t) < \rho + \rho'< 2\rho$$ and $u\circ\alpha_t$ is a $W^{1,2}$--curve it follows from \cite[Lemma 8.5]{LW15-Plateau} that there exists $w^+\in W^{1,2}(\Omega^+, X)$ with $$\Area(w^+)< 4C\cdot \rho^2$$ and such that $\trace(w^+)\circ \alpha_t^+ = u\circ\alpha_t$ and $\trace(w^+)|_{\tilde{\gamma}^+}$ is a constant speed parametrization of $\Gamma_j^-$. Moreover, since the continuous representative of $\trace(u)|_{\gamma^-\cup \alpha_t^-}$  satisfies $$\length(\trace(u)|_{\gamma^-\cup \alpha_t^-}) = \length(\Gamma_i^-) + \length(u\circ\alpha_t) < 2\rho$$ and is the trace of a Sobolev annulus it follows from the local quadratic isoperimetric inequality and from \cite[Lemma 4.8]{LW-intrinsic} and its proof that there exists $w^-\in W^{1,2}(\Omega^-, X)$ with $\trace(w^-) = \trace(u)|_{\gamma^-\cup \alpha_t^-}$ and such that $$\Area(w^-)< 4C\cdot \rho^2.$$

Again by \cite[Theorem 1.12.3]{KS93}, the map $v$ which coincides with $u$ on $M'$ and with $w^\pm$ on $\Omega^\pm$ belongs to $W^{1,2}(M^*, X)$ and satisfies 
\begin{equation}\label{eq:reduction-area-two}
\Area(v) = \Area(u) + \Area(w^-) + \Area(w^+) < \Area(u) + 8C\rho^2\leq \Area(u) + \eta.
\end{equation}
Moreover, on $\partial M^*\setminus \tilde{\gamma}^+$  the trace of $v$ coincides with $\trace(u)|_{\partial M \setminus \gamma^-}$;  on ${\tilde{\gamma}^+}$ it is a constant speed parametrization of $\Gamma_j^-$. In particular, $\trace(v)$ is a weakly monotone parametrization of $\Gamma$ and hence $v\in \Lambda(M^*,\Gamma, X)$. 

We can arrive at a contradiction as in the first case. Indeed, the Euler characteristic of $M^*$ satisfies $$\chi(M^*) = \chi(M') + 1 = \chi(M) + 2.$$ Thus, if $M^*$ is connected then it has genus $p-1$ and hence is a reduction of $M$. If $M^*$ has two connected components then it has genus $p$. So, if each component has non-empty boundary then $M^*$ is again a reduction of $M$. In both cases, inequality \eqref{eq:reduction-area-two} shows that $$a_p^*(\Gamma, X) \leq \Area(v)<\Area(u) + \eta,$$ which contradicts the assumption on $u$. Finally, if one of the components of $M^*$ has empty boundary then it has genus at least one. In this case the surface obtained by omitting this component is a reduction of $M$ and we obtain a contradiction in the same way as above. This concludes the proof of the proposition.
\end{proof}

\section{Proof of the main results}\label{sec:area-mins}

In this section, we prove the results stated in the introduction. We begin with Theorem~\ref{thm:Douglas-Plateau-intro}, let $(X,d)$ be a proper metric space admitting a local quadratic isoperimetric inequality and let $\Gamma$ be the disjoint union of $k\geq 1$ rectifiable Jordan curves in $X$. We furthermore let $M$ be a smooth compact and connected surface with $k$ boundary components and of genus $p\geq 0$. We first assume that $k+2p\geq 3$ and hence $M$ admits a hyperbolic metric. The Douglas condition \eqref{eq:Douglas-condition} implies that the family $\Lambda:=\Lambda(M, \Gamma, X)$ is not empty.

\bp\label{prop:good-limit}
 Let $(g_n)$ be a sequence of hyperbolic metrics on $M$. Suppose $(u_n)\subset \Lambda$ is a sequence which satisfies 
 \begin{equation}\label{eq:area-below-apstar}
  \sup_n \Area(u_n)< a_p^*(\Gamma, X)
 \end{equation}
and 
\begin{equation}\label{eq:energy-conv-to-some-m}
 \lim_{n\to\infty} E_+^2(u_n, g_n)= m
\end{equation}
 for some $m>0$.
Then there exist $v\in \Lambda$ and a hyperbolic metric $g$ on $M$ with the following property. After precomposing each $u_n$ with a diffeomorphism of $M$ and passing to a subsequence, the maps $u_n$ converge to $v$ in $L^2(M, X)$ and one has $$\lim_{n\to\infty} E_+^2(u_n, g) =m.$$
\ep

\begin{proof}
 Let $(g_n)$, $(u_n)$ and $m$ be as in the statement of the proposition. In view of \eqref{eq:area-below-apstar} and \eqref{eq:energy-conv-to-some-m}, 
Proposition~\ref{prop:bound-systole} implies that the relative systole of $(M, g_n)$ is bounded away from zero independently of $n$. By Theorem~\ref{thm:Mumford-cptness-boundary} there thus exist diffeomorphisms $\varphi_n\colon M\to M$ such that, after possibly passing to a subsequence, the Riemannian metrics $\varphi_n^* g_n$ converge smoothly to a hyperbolic metric $g$ on $M$.
The maps $v_n:= u_n\circ \varphi_n$ belong to $\Lambda$ and satisfy $$\lim_{n\to\infty} E_+^2(v_n, g) = m$$ because $\varphi_n$, when viewed as a map from $(M, g)$ to $(M, g_n)$, is $\lambda_n$--biLipschitz with $\lambda_n\to 1$ as $n$ tends to infinity.
It now follows from Lemma~\ref{lem:bound-L2-norm} and the metric-space valued version of the Rellich-Kondrachov compactness theorem (see \cite[Theorem 1.13]{KS93}) that there exists a subsequence $(v_{n_j})$ which converges in $L^2(M, X)$ to some map $v\in W^{1,2}(M, X)$. 

It remains to show that $v$ belongs to $\Lambda$. By Proposition~\ref{prop:equi-cont-traces} the sequence $(\trace(v_n))$ is equi-continuous. Therefore, by Arzel\`a-Ascoli theorem, a subsequence of $(\trace(v_n))$ converges uniformly to some continuous map $\gamma\colon \partial M\to X$. The map $\gamma$ is a weakly monotone parametrization of $\Gamma$ because it is the uniform limit of maps with this property. Finally, since the sequence $(\trace(v_n))$ converges in $L^2(\partial M, X)$ to $\trace(v)$ by \cite[Theorem 1.12.2]{KS93}, it follows that $\trace(v) = \gamma$ almost everywhere on $\partial M$. This shows that $v\in \Lambda$ and completes the proof.
\end{proof}

We can now finish the proof of our main theorem:

\begin{proof}[Proof of Theorem~\ref{thm:Douglas-Plateau-intro}]
We first assume that $k+2p\geq 3$ and claim that if $u\in \Lambda$ is such that $$\Area(u)< a_p^*(\Gamma, X)$$ then there exist $v\in\Lambda$ and a hyperbolic metric $g$ on $M$ such that $\Area(v)\leq \Area(u)$ and $v$ is infinitesimally isotropic with respect to $g$. Indeed, set $$\Lambda_u = \{ v\in \Lambda: \Area(v) \leq \Area(u)\}$$ and
$$m= \inf\{ E_+^2(v, g): \text{ $v\in \Lambda_u$, $g$ hyperbolic metric on $M$}\}$$
and choose a sequence $(u_n, g_n)$, where $u_n\in \Lambda_u$ and $g_n$ is a hyperbolic metric, such that $$\lim_{n\to\infty} E_+^2(u_n, g_n) = m.$$ 
By Proposition~\ref{prop:good-limit} there exist a map $v\in\Lambda$ and a hyperbolic metric $g$ on $M$ such that, after possibly precomposing each $u_n$ by a diffeomorphism of $M$ and passing to a subsequence, the maps $u_n$ converge to $v$ in $L^2(M,X)$ and $E_+^2(u_n, g) \to m$ as $n$ tends to infinity. By the lower semi-continuity of area and energy (see \cite[Corollaries 5.8 and 5.7]{LW15-Plateau}) we have that $\Area(v)\leq \Area(u)$ and  $E_+^2(v, g) \leq m$. This implies, in particular, that $v\in\Lambda_u$ and so $E_+^2(v, g) = m$. It now follows from the invariance of area under biLipschitz homeomorphisms and from Theorem~\ref{thm:energy-min-isotropic} that $v$ is infinitesimally isotropic with respect to $g$. This proves our claim. Notice that in \cite{LW15-Plateau} the lower semi-continuity results referred to above are proved for maps defined on open bounded subsets of $\R^2$. The corresponding results for maps defined on $M$ easily follow from this by decomposing $M$ into a disjoint union $M=U_1\cup\dots\cup U_L\cup N$, where $N$ is a set of measure zero and each $U_i$ is an open disc whose closure is contained in a conformal chart.

Now, let $(u_n)\subset\Lambda$ be an area minimizing sequence, thus $$\Area(u_n) \to a(M, \Gamma, X)=a_p(\Gamma, X)$$ as $n$ tends to infinity. Since the Douglas condition \eqref{eq:Douglas-condition} holds we may assume that $$\sup_n\Area(u_n)< a_p^*(\Gamma, X).$$
By the claim above there exist a sequence $(v_n)\subset\Lambda$ and a sequence $(g_n)$ of hyperbolic metrics on $M$ such that $\Area(v_n)\leq \Area(u_n)$ and $v_n$ is infinitesimally isotropic with respect to $g_n$ for every $n$. In particular, $(v_n)$ is an area minimizing sequence and the sequence of energies $E_+^2(v_n, g_n)$ is uniformly bounded by \eqref{eq:energy-bdd-by-area-inf-isotropic}. Proposition~\ref{prop:good-limit} shows that there exists a map $v\in\Lambda$ such that, after possibly precomposing each $v_n$ by a diffeomorphism and passing to a subsequence, the maps $v_n$ converge to $v$ in $L^2(M, X)$ and that the energies $E_+^2(v_n, g)$ are uniformly bounded for some (and thus every fixed) hyperbolic metric $g$ on $M$. By the lower semi-continuity of area  (see \cite[Corollary 5.8]{LW15-Plateau}) we have $$\Area(v)\leq \liminf_{n\to\infty}\Area(v_n) = a(M, \Gamma, X)$$ and hence $\Area(v) = a(M, \Gamma, X)$. This shows the existence of an area minimizer in $\Lambda$. 
Applying the claim at the beginning of the proof again, we obtain the existence of a hyperbolic metric $g'$ and an area minimizer $v'$ in $\Lambda$ such that $v'$ is infinitesimally isotropic with respect to $g'$. This concludes the proof of Theorem~\ref{thm:Douglas-Plateau-intro} in the case $k+2p\geq 3$.

The case $k+2p= 2$ works analogously. Then $M$ is diffeomorphic to a cylinder and we work with flat metrics, thus Riemannian metrics $g$ for which $(M, g)$ has constant curvature $0$ and $\partial M$ is geodesic. We furthermore normalize so that $(M, g)$ has area equal to $1$. Theorem~\ref{thm:Mumford-cptness-boundary} has a natural analogue in this case (see \cite[Theorem 4.4.1]{DHT10} for the version for closed surfaces). The analogue of Proposition~\ref{prop:bound-systole} for the case of flat metrics is proved similarly and relies on the existence of a suitable flat collar which is elementary in this case. The rest of the proof of the theorem in the case $k+2p=2$ remains unchanged.

Finally, we note that the remaining case $k+2p = 1$ is exactly the classical problem of Plateau treated in \cite{LW15-Plateau}.
\end{proof}

\begin{proof}[Proof of Theorem~\ref{thm:area-min-bounded-genus}]
 Let $X$, $\Gamma$, $k$, $p$ be as in the statement of the theorem and assume that $\hat{a}_p(\Gamma, X)<\infty$. Up to diffeomorphisms, the family $\hat{\mathcal{M}}(k,p)$ contains only finitely many different smooth surfaces. We denote these by $M_1,\dots, M_N$ and thus have $$\hat{a}_p(\Gamma, X) = \min\{a(M_j, \Gamma, X): j=1,\dots, N\}.$$ Among those $M_j$ with $a(M_j, \Gamma, X) = \hat{a}_p(\Gamma, X)$ choose one for which $M_j$ has the largest Euler characteristic and denote it by $M$. Let $(u_n)\subset \Lambda(M, \Gamma, X)$ be an area minimizing sequence, thus 
 \begin{equation}\label{eq:area-min-bounded-genus}
 \Area(u_n) \to a(M, \Gamma, X) = \hat{a}_p(\Gamma, X)
 \end{equation}
  as $n$ tends to infinity.
 
Let $M^1,\dots, M^m$ be the connected components of $M$. There exist a subsequence $(u_{n_l})$ and a partition $\Gamma = \Gamma^1\cup\dots\cup\Gamma^m$ of $\Gamma$ into unions of Jordan curves such that $\trace(u_{n_l})(\partial M^i) = \Gamma^i$ for all $l\in\N$ and all $i$. Thus $u_{n_l}|_{M^i} \in\Lambda(M^i, \Gamma^i, X)$ and it follows with \eqref{eq:area-min-bounded-genus} that $$\hat{a}_p(\Gamma, X) = \sum_{i=1}^m a(M^i, \Gamma^i, X).$$
Let $p_i$ be the genus of $M^i$. We claim that 
\begin{equation*}
a_{p_i}(\Gamma^i, X) = a(M^i, \Gamma^i, X) < a_{p_i}^*(\Gamma^i, X)
\end{equation*}
 for all $i$. We argue by contradiction and assume that this is wrong for some $i$. There thus exists a reduction $M^i_*$ of $M^i$ such that $$a(M^i_*, \Gamma^i, X) \leq a(M^i, \Gamma^i, X),$$ and in fact equality holds. Let $M^*$ be the surface obtained from $M$ by replacing the connected component $M^i$ by $M^i_*$. It is clear that $M^*\in \hat{\mathcal{M}}(k,p)$. Moreover, we have $a(M^*, \Gamma, X) = \hat{a}_p(\Gamma, X)$ but $M^*$ has strictly larger Euler characteristic than $M$, which contradicts the choice of $M$. This proves our claim.

By the claim and Theorem~\ref{thm:Douglas-Plateau-intro} there exist for each $i$ an element $u^i\in\Lambda(M^i, \Gamma^i, X)$ and a Riemannian metric $g_i$ on $M^i$ such that $$\Area(u^i) = a(M^i, \Gamma^i, X)$$ and $u^i$ is infinitesimally isotropic with respect to $g_i$. Then the map $u$ coinciding with $u^i$ on $M^i$ belongs to $\Lambda(M, \Gamma, X)$ and satisfies $\Area(u) = \hat{a}_p(\Gamma, X)$. Moreover, $u$ is infinitesimally isotropic with respect to the Riemannian metric $g$ on $M$ which agrees with $g_i$ on each $M^i$. This completes the proof.
 \end{proof}

We finally indicate how to deduce Theorem~\ref{thm:regularity-sol-Douglas-Plateau} from the results in \cite{LW15-Plateau}.

\begin{proof}[Proof of Theorem~\ref{thm:regularity-sol-Douglas-Plateau}]
 Let $(U,\psi)$ be a conformal chart in the interior of $M$ satisfying $\psi(U) = D$. For almost every $r\in(0,1)$ the map $v(z):= u\circ\psi^{-1}(rz)$ belongs to $W^{1,2}(D, X)$ and satisfies $$\Area(v) = \inf\{\Area(w): \text{$w\in W^{1,2}(D, X)$, $\trace(w) = \trace(v)$}\}.$$ Moreover, $v$ is infinitesimally isotropic with respect to $g_{\rm Eucl}$ and hence $\sqrt{2}$--quasi\-conformal in the sense of \cite{LW15-Plateau}. It thus follows from \cite[Theorem 8.2]{LW15-Plateau} that $v\in W^{1,q}_{\rm loc}(D, X)$ for some $q>2$ and that the continuous representative $\bar{v}$ of $v$ is locally $\alpha$--H\"older continuous with $\alpha = (8\pi C)^{-1}$. In particular, $\bar{v}$ satisfies Lusin's property (N). The value of $q$ only depends on $C$, see the proof of \cite[Theorem 8.2]{LW15-Plateau}. Since the map $z\mapsto \psi^{-1}(rz)$ is biLipschitz for fixed $r\in(0,1)$ this establishes statement (i) and the first part of statement (ii) of our theorem. In what follows, we denote by $\bar{u}$ the representative of $u$ which is continuous in the interior of $M$.
 
In order to prove the second part of statement (ii) let $(U, \psi)$ be a conformal chart around a boundary point of $M$ with image $$\psi(U) = \{z\in D: \operatorname{Im}(z)\geq 0\}.$$ For $r\in(0,1)$ define $$D_r^+:= \{z\in \C: \text{$|z|<r$ and $\operatorname{Im}(z)>0$}\}.$$ Then for almost every $r$ the map $v\colon D_r^+\to X$ given by $v:= \bar{u}\circ\psi^{-1}|_{D_r^+}$ is in $W^{1,2}(D_r^+, X)$, is infinitesimally isotropic with respect to $g_{\rm Eucl}$ and satisfies $$\Area(v) = \inf\{\Area(w): \text{$w\in W^{1,2}(D_r^+, X)$, $\trace(w) = \trace(v)$}\},$$ and $\trace(v)$ has a continuous representative. It thus follows from \cite[Theorem 9.1]{LW15-Plateau} that $v$ extends continuously to the boundary of $D_r^+$. This shows that $\bar{u}$ has a continuous extension to all of $\partial M$, thus proving the second part of statement (ii).

Finally, statement (iii) follows almost as in the proof of \cite[Theorem 9.3]{LW15-Plateau}. Notice that the $3$--point condition assumed in that proof is not needed provided the value of $r_0>0$ appearing therein is chosen sufficiently small.
\end{proof}

\section{Courant's condition of cohesion}\label{sec:cohesion}

We recall the condition of cohesion introduced by Courant \cite{Cou40}  and used in \cite{Shi39} and \cite{TT88}. This condition is for example satisfied when the maps are incompressible in the sense of Schoen-Yau \cite{SY79}. We then prove the existence of energy minimizers in proper metric spaces under the condition of cohesion.

Let $X$ be a complete metric space and $M$ a smooth compact and connected surface.

\bd
A map $u\colon M\to X$ is called $\eta$--cohesive, $\eta>0$, if $u$ is continuous and $$\length(u\circ c)\geq \eta$$ for every non-contractible closed curve $c$ in $M$.
A family $\mathcal{F}$ of maps from $M$ to $X$ is said to satisfy the condition of cohesion if there exists $\eta>0$ such that each $u\in\mathcal{F}$ is $\eta$--cohesive. 
\ed

Let $\Gamma$ be the disjoint union of $k\geq 1$ rectifiable Jordan curves in $X$ and suppose $M$ has $k$ boundary components. Set $$e(M,\Gamma, X):= \inf\{E_+^2(u,g): \text{ $u\in\Lambda(M,\Gamma, X)$, $g$ Riemannian metric on $M$}\}.$$
An {\it energy minimizing sequence} in $\Lambda(M,\Gamma, X)$ is a sequence of pairs $(u_n, g_n)$ of maps $u_n\in\Lambda(M,\Gamma, X)$ and Riemannian metrics $g_n$ on $M$ satisfying $$E_+^2(u_n, g_n) \to e(M,\Gamma, X)$$ as $n$ tends to infinity. 

The following theorem generalizes \cite{Shi39}, \cite{Cou40} and \cite{TT88} to the setting of proper metric spaces.

\bt\label{thm:cond-cohesion-energy-min}
 Let $X$ be a proper metric space and let $\Gamma$ be the disjoint union of $k\geq 1$ rectifiable Jordan curves in $X$. Let $M$ be a smooth compact and connected surface with $k$ boundary components. If there is an energy minimizing sequence in $\Lambda(M, \Gamma, X)$ satisfying the condition of cohesion then there exist $u\in\Lambda(M, \Gamma, X)$ and a Riemannian metric $g$ on $M$ such that $$E_+^2(u,g)=e(M,\Gamma, X).$$ For any such $u$ and $g$ the map $u$ is infinitesimally isotropic with respect to $g$.
\et

The Riemannian metric $g$ can be chosen in such a way that $(M,g)$ has constant curvature $-1$, $0$, $1$ and that $\partial M$ is geodesic.

In the generality of metric spaces, energy minimizers in $\Lambda(M,\Gamma, X)$ with respect to the Reshetnyak energy $E_+^2$ need not be minimizers of the parametrized Hausdorff area, see \cite[Proposition 11.6]{LW15-Plateau}. However, one can show that they are minimizers of the parametrized area coming from the so-called inscribed Riemannian area, see \cite{FW-Morrey} and compare with \cite{LW17-en-area} where this is proved when $M$ is the disc. In particular, if $X$ also admits a local quadratic isoperimetric inequality then maps $u$ as in Theorem~\ref{thm:cond-cohesion-energy-min} are locally H\"older on $M\setminus \partial M$ and extend continuously to $\partial M$ by the analogue of Theorem~\ref{thm:regularity-sol-Douglas-Plateau} for the inscribed Riemannian area. Notice that if $X$ has the so-called property (ET) introduced in \cite{LW15-Plateau} then the inscribed Riemannian area agrees with the parametrized Hausdorff area. 

We finally mention that one can combine the arguments from the proof of Theorem~\ref{thm:cond-cohesion-energy-min} with a metric space version of the Morrey $\varepsilon$--conformality lemma to obtain the existence of a Hausdorff area minimizer under the condition of cohesion for an area minimizing sequence, see \cite{FW-Morrey}. However, unlike in Theorem~\ref{thm:Douglas-Plateau-intro}, we do not know how to obtain a good parametrization for such an area minimizer.

We turn to the proof of the theorem. The following provide analogues of the propositions proved in Sections~\ref{sec:equi-cont-traces} and \ref{sec:lower-bound-rel-sys}.

\bp\label{prop:equi-cont-cohesive}
 Let $g$ be a Riemannian metric on $M$. Then for every $\eta>0$ and $K>0$ the family $$\{u|_{\partial M}: \text{ $u\in\Lambda(M, \Gamma, X)$ is $\eta$--cohesive and $E_+^2(u,g)\leq K$}\}$$ is equi-continuous.
\ep

\begin{proof}
 This follows from the same arguments as those used at the beginning of the proof of Proposition~\ref{prop:equi-cont-traces}. 
 The value of $\rho$ in that proof is replaced by $\rho:= \min\{\varepsilon, \frac{\eta}{2}\}$ and one uses the fact that $u$ is $\eta$--cohesive to show that the curve $\gamma^-$ in the proof must satisfy $\length(u\circ\gamma^-)\leq \rho\leq \varepsilon$. Indeed, otherwise the curve $c$ obtained by concatenating the curves $\beta_r$ and $\gamma^+$ (appearing in the proof) provides a non-contractible curve in $M$ such that $$\length(u\circ c) = \length(u\circ\beta_r) + \length(u\circ\gamma^+)<2\rho\leq \eta.$$ This contradicts the assumption that $u$ is $\eta$--cohesive.
\end{proof}

\bp\label{prop:lower-bound-rel-sys-cohesive}
 For every $\eta>0$ and $K>0$ there exists $\varepsilon>0$ with the following property. If $u\in\Lambda(M, \Gamma, X)$ is $\eta$--cohesive and $g$ is a hyperbolic metric on $M$ such that $E_+^2(u,g)\leq K$ then the relative systole of $(M,g)$ is bounded from below by $\varepsilon$.
\ep

\begin{proof}
 This uses the same arguments as those in the proof of Proposition~\ref{prop:bound-systole}. The value of $\rho$ appearing therein is replaced by $\rho:=\frac{\eta}{2}$. If the relative systole of $(M,g)$ is smaller than the $\varepsilon$ in that proof then one obtains a contradiction with the fact that $u$ is $\eta$--cohesive. Indeed, let $\alpha_t$ be the curve in that proof. If $\alpha_t$ is a closed curve then \eqref{eq:length-short-parallel} already yields a contradiction. If $\alpha_t$ has endpoints on $\partial M$ then the concatenation $c$ of $\alpha_t$ with the curve $\gamma^-$ appearing in the proof is a non-contractible curve such that $$\length(u\circ c) = \length(u\circ\gamma^-) + \length(u\circ\alpha_t)<2\rho\leq \eta,$$ which is again a contradiction.
\end{proof}

\begin{proof}[Proof of Theorem~\ref{thm:cond-cohesion-energy-min}]
 We only sketch the proof for the case that $M$ has strictly negative Euler characteristic, the case of a cylinder being analogous and the case of a disc appearing in \cite[Theorem 7.6]{LW15-Plateau}.
 
Let $(u_n, g_n)$ be an energy minimizing sequence in $\Lambda(M,\Gamma, X)$ which satisfies the condition of cohesion for some $\eta>0$. There exists a hyperbolic metric on $M$ which is conformally equivalent to $g_n$. By the conformal invariance of the Reshetnyak energy, we may thus assume that each $g_n$ is hyperbolic. By Proposition~\ref{prop:lower-bound-rel-sys-cohesive} the relative systole of $(M,g_n)$ is bounded away from zero independently of $n$.
 
Using Proposition~\ref{prop:equi-cont-cohesive} instead of Proposition~\ref{prop:equi-cont-traces} and arguing exactly as in the proof of Proposition~\ref{prop:good-limit} one obtains: there exist $v\in\Lambda(M, \Gamma, X)$ and  a hyperbolic metric $g$ on $M$ with the following property. After precomposing each $u_n$ with a suitable diffeomorphism of $M$ and passing to a subsequence, the maps $u_n$ converge to $v$ in $L^2(M, X)$ and $$\lim_{n\to\infty} E_+^2(u_n, g)= e(M,\Gamma, X).$$ By the lower semi-continuity of energy we have $E_+^2(v,g)\leq e(M,\Gamma, X)$ and thus equality holds. This shows the existence of an energy minimizing pair $(v,g)$.

Finally, Theorem~\ref{thm:energy-min-isotropic} shows that for any $u\in\Lambda(M,\Gamma, X)$ and for any Riemannian metric $g$ with $E_+^2(u,g) = e(M,\Gamma, X)$ the map $u$ is infinitesimally isotropic with respect to $g$.
\end{proof}

\def\cprime{$'$} \def\cprime{$'$} \def\cprime{$'$}


\begin{thebibliography}{10}

\bibitem{AS60}
Lars~V. Ahlfors and Leo Sario.
\newblock {\em Riemann surfaces}.
\newblock Princeton Mathematical Series, No. 26. Princeton University Press,
  Princeton, N.J., 1960.

\bibitem{Amb90}
Luigi Ambrosio.
\newblock Metric space valued functions of bounded variation.
\newblock {\em Ann. Scuola Norm. Sup. Pisa Cl. Sci. (4)}, 17(3):439--478, 1990.

\bibitem{AK00}
Luigi Ambrosio and Bernd Kirchheim.
\newblock Currents in metric spaces.
\newblock {\em Acta Math.}, 185(1):1--80, 2000.

\bibitem{Bal97}
Keith Ball.
\newblock An elementary introduction to modern convex geometry.
\newblock In {\em Flavors of geometry}, volume~31 of {\em Math. Sci. Res. Inst.
  Publ.}, pages 1--58. Cambridge Univ. Press, Cambridge, 1997.

\bibitem{Bus10}
Peter Buser.
\newblock {\em Geometry and spectra of compact {R}iemann surfaces}.
\newblock Modern Birkh\"{a}user Classics. Birkh\"{a}user Boston, Inc., Boston,
  MA, 2010.
\newblock Reprint of the 1992 edition.

\bibitem{Cou40}
R.~Courant.
\newblock The existence of minimal surfaces of given topological structure
  under prescribed boundary conditions.
\newblock {\em Acta Math.}, 72:51--98, 1940.

\bibitem{DHT10}
Ulrich Dierkes, Stefan Hildebrandt, and Anthony~J. Tromba.
\newblock {\em Global analysis of minimal surfaces}, volume 341 of {\em
  Grundlehren der Mathematischen Wissenschaften [Fundamental Principles of
  Mathematical Sciences]}.
\newblock Springer, Heidelberg, second edition, 2010.

\bibitem{Dou31}
Jesse Douglas.
\newblock Solution of the problem of {P}lateau.
\newblock {\em Trans. Amer. Math. Soc.}, 33(1):263--321, 1931.

\bibitem{Dou39}
Jesse Douglas.
\newblock Minimal surfaces of higher topological structure.
\newblock {\em Ann. of Math. (2)}, 40(1):205--298, 1939.

\bibitem{EG92}
Lawrence~C. Evans and Ronald~F. Gariepy.
\newblock {\em Measure theory and fine properties of functions}.
\newblock Studies in Advanced Mathematics. CRC Press, Boca Raton, FL, 1992.

\bibitem{FW-Morrey}
Martin Fitzi and Stefan Wenger.
\newblock Morrey $\varepsilon$-conformality in metric spaces and applications.
\newblock {\em in preparation}.

\bibitem{Haj96}
Piotr Haj{\l}asz.
\newblock Sobolev spaces on an arbitrary metric space.
\newblock {\em Potential Anal.}, 5(4):403--415, 1996.

\bibitem{Heb99}
Emmanuel Hebey.
\newblock {\em Nonlinear analysis on manifolds: {S}obolev spaces and
  inequalities}, volume~5 of {\em Courant Lecture Notes in Mathematics}.
\newblock New York University, Courant Institute of Mathematical Sciences, New
  York; American Mathematical Society, Providence, RI, 1999.

\bibitem{HKST15}
Juha Heinonen, Pekka Koskela, Nageswari Shanmugalingam, and Jeremy Tyson.
\newblock {\em Sobolev spaces on metric measure spaces}, volume~27 of {\em New
  Mathematical Monographs}.
\newblock Cambridge University Press, Cambridge, 2015.

\bibitem{Jos85}
J{\"u}rgen Jost.
\newblock Conformal mappings and the {P}lateau-{D}ouglas problem in
  {R}iemannian manifolds.
\newblock {\em J. Reine Angew. Math.}, 359:37--54, 1985.

\bibitem{Jos97}
J\"{u}rgen Jost.
\newblock Generalized {D}irichlet forms and harmonic maps.
\newblock {\em Calc. Var. Partial Differential Equations}, 5(1):1--19, 1997.

\bibitem{Jos06-cptRiem}
J\"{u}rgen Jost.
\newblock {\em Compact {R}iemann surfaces}.
\newblock Universitext. Springer-Verlag, Berlin, third edition, 2006.
\newblock An introduction to contemporary mathematics.

\bibitem{Kar07}
M.~B. Karmanova.
\newblock Area and co-area formulas for mappings of the {S}obolev classes with
  values in a metric space.
\newblock {\em Sibirsk. Mat. Zh.}, 48(4):778--788, 2007.

\bibitem{Kir94}
Bernd Kirchheim.
\newblock Rectifiable metric spaces: local structure and regularity of the
  {H}ausdorff measure.
\newblock {\em Proc. Amer. Math. Soc.}, 121(1):113--123, 1994.

\bibitem{KS93}
Nicholas~J. Korevaar and Richard~M. Schoen.
\newblock Sobolev spaces and harmonic maps for metric space targets.
\newblock {\em Comm. Anal. Geom.}, 1(3-4):561--659, 1993.

\bibitem{LW15-Plateau}
Alexander Lytchak and Stefan Wenger.
\newblock Area minimizing discs in metric spaces.
\newblock {\em Arch. Ration. Mech. Anal.}, 223(3):1123--1182, 2017.

\bibitem{LW17-en-area}
Alexander Lytchak and Stefan Wenger.
\newblock Energy and area minimizers in metric spaces.
\newblock {\em Adv. Calc. Var.}, 10(4):407--421, 2017.

\bibitem{LW-intrinsic}
Alexander Lytchak and Stefan Wenger.
\newblock Intrinsic structure of minimal discs in metric spaces.
\newblock {\em Geom. Topol.}, 22(1):591--644, 2018.

\bibitem{MZ10}
Chikako Mese and Patrick~R. Zulkowski.
\newblock The {P}lateau problem in {A}lexandrov spaces.
\newblock {\em J. Differential Geom.}, 85(2):315--356, 2010.

\bibitem{Mor48}
Charles~B. Morrey, Jr.
\newblock The problem of {P}lateau on a {R}iemannian manifold.
\newblock {\em Ann. of Math. (2)}, 49:807--851, 1948.

\bibitem{Nik79}
I.~G. Nikolaev.
\newblock Solution of the {P}lateau problem in spaces of curvature at most
  {$K$}.
\newblock {\em Sibirsk. Mat. Zh.}, 20(2):345--353, 459, 1979.

\bibitem{OvdM14}
Patrick Overath and Heiko von~der Mosel.
\newblock Plateau's problem in {F}insler 3-space.
\newblock {\em Manuscripta Math.}, 143(3-4):273--316, 2014.

\bibitem{Rad30}
Tibor Rad{\'o}.
\newblock On {P}lateau's problem.
\newblock {\em Ann. of Math. (2)}, 31(3):457--469, 1930.

\bibitem{Res97}
Yu.~G. Reshetnyak.
\newblock Sobolev classes of functions with values in a metric space.
\newblock {\em Sibirsk. Mat. Zh.}, 38(3):657--675, iii--iv, 1997.

\bibitem{Res06}
Yu.~G. Reshetnyak.
\newblock On the theory of {S}obolev classes of functions with values in a
  metric space.
\newblock {\em Sibirsk. Mat. Zh.}, 47(1):146--168, 2006.

\bibitem{SY79}
R.~Schoen and Shing~Tung Yau.
\newblock Existence of incompressible minimal surfaces and the topology of
  three-dimensional manifolds with nonnegative scalar curvature.
\newblock {\em Ann. of Math. (2)}, 110(1):127--142, 1979.

\bibitem{Shi39}
Max Shiffman.
\newblock The {P}lateau problem for minimal surfaces of arbitrary topological
  structure.
\newblock {\em Amer. J. Math.}, 61:853--882, 1939.

\bibitem{TT88}
Friedrich Tomi and Anthony~J. Tromba.
\newblock Existence theorems for minimal surfaces of nonzero genus spanning a
  contour.
\newblock {\em Mem. Amer. Math. Soc.}, 71(382):iv+83, 1988.

\end{thebibliography}
\end{document}